\documentclass[onecolumn]{autart}
\usepackage{}    
\usepackage[colorlinks,linkcolor=black,anchorcolor=red,citecolor=red]{hyperref}%

\usepackage{graphicx}           
\usepackage{dcolumn}            
\usepackage{bm}                 
\usepackage{bbm}
\usepackage{graphics}           
\usepackage{epsfig}             
\usepackage{times}              
\usepackage[fleqn]{amsmath}
\usepackage{amssymb}
\usepackage{extarrows}
\usepackage{amsfonts}
\usepackage{mathrsfs}
\usepackage{fancyhdr}
\usepackage{color}
\usepackage{subfigure}
\usepackage{enumerate}
\usepackage [latin1]{inputenc} %
\allowdisplaybreaks[4] 

\newtheorem{theorem}{Theorem}[section] %
\newtheorem{lemma}{Lemma}[section]

\newtheorem{proposition}[theorem]{Proposition}

\newtheorem{definition}{Definition}[section]

\newtheorem{remark}{Remark}[section]

\newcommand{\diff}[1]{\operatorname{d}\!{#1}} 
\newcommand{\divi}{\operatorname{div}\!}

\begin{document}

\begin{frontmatter}

\title{Approximations of Lyapunov functions for ISS analysis of a class of nonlinear parabolic PDEs}

\author{Jun Zheng$^{1}$}\ead{zhengjun2014@aliyun.com},
\author{Guchuan Zhu$^{2}$}\ead{guchuan.zhu@polymtl.ca}

\address{$^{1}${School of Mathematics, Southwest Jiaotong University,
        Chengdu 611756, Sichuan, China}\\
        $^{2}$Department of Electrical Engineering, Polytechnique Montr\'{e}al, P.O. Box 6079, Station Centre-Ville, Montreal, QC, Canada H3T 1J4}

\begin{abstract}
{This paper addresses the input-to-state stability (ISS) and integral input-to-state stability (iISS) for a class of nonlinear higher dimensional parabolic partial differential equations (PDEs) with different types of boundary disturbances (Robin or Neumann or Dirichlet) from different spaces by means of approximations of Lyapunov functions. Specifically, by constructing approximations of (coercive and non-coercive) ISS Lyapunov functions we establish: (i) the ISS and iISS in $L^1$-norm (and weighted $L^1$-norm) for PDEs with boundary disturbances from $L^q_{loc}(\mathbb{R}_+;L^1(\partial\Omega))$-space for any $q\in [1,+\infty]$; (ii) the iISS in $L^1$-norm (and weighted $L^1$-norm) for PDEs with boundary disturbances from $L^\Phi_{loc}(\mathbb{R}_+;L^1(\partial\Omega))$-space for certain Young function $\Phi$; and (iii) the ISS and iISS in $L^\Phi$-norm (and weighted $K_\Phi$-class) for PDEs with boundary disturbances from $L^q_{loc}(\mathbb{R}_+;K_\Phi(\partial\Omega))$-class for any $q\in [1,+\infty]$  and certain Young function  $\Phi$. {The ISS properties stated} in (ii) and (iii) are assessed in the framework of Orlicz space or Orlicz class.}
\end{abstract}
\end{frontmatter}
\tableofcontents 
\section{Introduction}\label{Sec: Introduction}
In the past few years, the input-to-state stability (ISS) theory for infinite dimensional systems governed by partial differential equations (PDEs) has drawn much attention in the literature of PDE control (see, e.g., \cite{Mironchenko:2019b} for a comprehensive survey).
{From the point of view of regularity theory of PDEs, ISS analysis for specific PDEs amounts essentially to establishing \emph{a priori} ISS estimates of the solutions with respect to (w.r.t.) disturbances.} Compared to the classical regularity estimates of solutions to PDEs with inhomogeneous boundary conditions (see, e.g., \cite{Evans:2010,Ladyzhenskaya:1968,Wu2006}), ISS estimates for PDEs with boundary disturbances often have much more specific forms as indicated in \cite{Zheng:2019b}, {because they must reflect some particular characteristics of the solutions in the presence of external disturbances.} Nevertheless, it is expected that, with more technical treatments, {the methods developed in classical regularity theory} can be applied to ISS analysis for PDEs with boundary disturbances. {In line with this perspective, we can find that:} (i) the monotonicity-based method has been used for ISS of monotone parabolic systems with Dirichlet boundary disturbances \cite{Mironchenko:2019}, (ii) the De~Giorgi iteration has been used for ISS of parabolic PDEs with Dirichlet boundary disturbances \cite{Zheng:201702}, (iii) the variations of Sobolev embedding inequalities were used for ISS of parabolic PDEs with Robin or Neumann boundary disturbances \cite{Zheng:201804,Zheng:201802}, and (iv) the maximum principle was used for ISS of parabolic equations with different types of boundary disturbances satisfying compatibility conditions \cite{Zheng:2019b,Zheng:2020CDC}. It should be mentioned that {semigroups and functional analytic methods} are effective for ISS analysis of linear or certain nonlinear PDE systems \cite{Jacob:2018CDC,Jacob:2019,Jacob:2018_SIAM,Jacob:2018_JDE,Schwenninger:2019}, and the spectral approach or finite-difference schemes can be used for ISS analysis of linear PDEs  \cite{Karafyllis:2014,Karafyllis:2016a,karafyllis2017siam,Karafyllis:2017,Lhachemi:201901,Lhachemi:201902}. Besides, {other methods in classical regularity theory, such as Morser iteration, Rothe's method, Galerkin's method, etc., may also be adopted} for the establishment of local or global ISS estimates for specific nonlinear PDEs with boundary disturbances.

It should be noticed that combined with different schemes and techniques as aforementioned, the Lyapunov {method is widely used in ISS analysis for PDE systems, such as parabolic systems with boundary disturbances \cite{Mironchenko:2019,Zheng:201804,Zheng:201802,Zheng:2019b,Zheng:2020CDC}, semilinear parabolic PDEs and diagonal parabolic systems with Neumann boundary disturbances \cite{Schwenninger:2019}, etc.} The ISS of certain first order hyperbolic systems was also established by the Lyapunov method in \cite{Tanwani:2017}.
It has been shown that one can obtain the ISS in $L^{q_1}$-norm (${q_1}\in [2,+\infty)$) for PDEs with boundary disturbances from $L^{q_2}$-space (${q_2}\in [2,+\infty]$) (see, e.g., \cite{Jacob:2019,Mironchenko:2019,Schwenninger:2019,Zheng:201804,Zheng:201802}) by the Lyapunov arguments combined with other methods. However, it is still challenging to establish the ISS in $L^{q_1}$-norm for PDEs with boundary (or in-domain) disturbances from $L^{q_2}$-space with any ${q_1}\in [1,+\infty]$ and any ${q_2}\in [1,+\infty]$ {by using solely the Lyapunov arguments. Particularly, to the best of the authors' knowledge, there are no results reported in the literature on establishing} the ISS estimates in $L^1$-norm for nonlinear PDEs with boundary (or in-domain) disturbances from $L^1$-space using the Lyapunov method. {This motivates the work presented in this paper.}


Specifically, {the first aim} of this paper is to establish ISS estimates in $L^{q_1}$-norm for PDEs with boundary (or in-domain) disturbances from $L^{q_2}$-space with any ${q_1}\in [1,+\infty]$ and any ${q_2}\in [1,+\infty]$. {Note that for} $q_1\geq 2$, a method based on non-coercive ISS Lyapunov functions was proposed  in \cite{Jacob:2019,Mironchenko:2018} to deal with different types of boundary disturbances, particularly, for linear parabolic PDEs with Dirichlet boundary disturbances.
 {In addition, the ISS in weighted $L^{q_1}$-norm (${q_1}=2$ or ${q_1}=+\infty$) (or weighted $L^1$-norm) for linear PDEs governed by Sturm-Liouville operators with $L^\infty$-disturbances (or $L^\infty$-Dirichlet disturbances) were established in \cite{karafyllis2017siam} by the eigenfunction expansion and the finite difference scheme.} However, constructing appropriate Lyapunov functions for the establishment of ISS estimates in $L^{q_1}$-norm for PDEs with boundary (or in-domain) disturbances from $L^{q_2}$-space may encounter serious difficulties when $q_1,q_2\in [1,2)$. The first obstacle is due to the well-posedness in the case {where} the inputs belong only to, e.g., $L^1$-space, since it is an not easy task to prove the existence of a ``solution'' (even in a weak sense) for PDEs (see remarks in Section~\ref{remarks existence}). The second obstacle lies at the construction of a Lyapunov function in the case of $q_1\in [1,2)$. Usually, when $q_1\geq 2$, a coercive Lyapunov candidate can be chosen as $V(u)=\|u\|^{q_1}_{L^{q_1}}$ for, e.g., parabolic PDEs with different boundary disturbances,  where $u$ is the solution of the considered PDEs (see, e.g., \cite{Jacob:2018_SIAM,Mironchenko:2019,Schwenninger:2019,Zheng:201702,Zheng:201804,Zheng:201802}) and a non-coercive ISS-Lyapunov function can be constructed for linear parabolic PDEs with Dirichlet boundary disturbances (see, e.g., \cite{Jacob:2019,Mironchenko:2018}). When  $q_1\in [1,2)$, the choice of $V(u)=\|u\|^{q_1}_{L^{q_1}}$ will fail to be effective for ISS analysis. 



We note that to overcome the first obstacle, the notion of ``entropy solution'' or ``renormalized solution'' can be used in well-posedness assessment for PDEs with inputs from $L^{1}$-space,{which is a topic we will consider in our future work}. In this paper, we focus on overcoming the second obstacle, i.e., we intend to establish ISS and iISS estimates in different norms for a class of nonlinear parabolic equations with boundary disturbances from different spaces by the Lyapunov method. More precisely, {we construct approximations of the desired Lyapunov functions that can be used to establish:} (i) the $L^q$-ISS\footnote{The specific definitions of $L^q$-ISS and $L^\Phi$-ISS are provided in Section \ref{Problem setting}.} ($q\in [1,+\infty]$) in $L^1$-norm (or weighted $L^1$-norm), (ii) the $L^q$-ISS ($q\in [1,+\infty]$) in $L^\Phi$-norm or $E_\Phi$-class (or, weighted $L^\Phi$-norm or weighted $E_\Phi$-class), and (iii) the $L^\Phi$-ISS in $L^1$-norm (or weighted $L^1$-norm), for a class of nonlinear PDEs with different boundary disturbances over higher dimensional domains, where $\Phi$ is a Young function. {}{In particular, as the second aim of this paper,} the results of (ii) and (iii) {}{are established} in the setting of Orlicz space, which is a generalization of $L^q$-space. The Orlicz space (or Orlicz class) has significant applications in mathematics and physics, particularly, in fluid mechanics. For example, it is well-known that the Lebesgue and Sobolev spaces can provide an appropriate functional framework for the study of existence and regularity of solutions to the nonlinear PDE model of non-Newtonian fluids when the nonlinear function is of power type as in \cite{Bird:1987}, while the Orlicz space is well suited for systems with nonlinearities having a non-polynomial type as in the Eyring-Prandtl model \cite{Eyring:1936}. The ISS in the framework of Orlicz spaces was first conducted in \cite{Jacob:2018_SIAM} for infinite dimensional systems with external disturbances, which can be seen as a generalization of classical stability notions in the control theory. Particularly, several characterizations of ISS for linear infinite dimensional systems with disturbances from $L^{q_2}$-space ($q_2\in[1,+\infty]$) and Orlicz space have been provided  in \cite{Jacob:2018_SIAM}, which also motivates us to address the problem in the setting of Orlicz spaces.

The main contribution of this paper includes:
\begin{enumerate}[(i)]
\item the construction of approximations for coercive or non-coercive Lyapunov functions that can be used to conduct ISS and iISS analysis in different norms for nonlinear PDE systems with boundary disturbances;
\item the establishment of $L^q$-ISS and $L^\Phi$-ISS in $L^1$-norm (or weighted $L^1$-norm) for nonlinear PDEs with boundary disturbances from $L^q$-space {}{and} Orlicz space, {}{respectively}, where $q\in [1,+\infty]$ and $\Phi$ is a certain Young function;
\item the establishment of $L^q$-ISS in $L^\Phi$-norm (or weighted $K_\Phi$-class) for nonlinear PDEs with boundary disturbances {}{from $L^q$-space}, where $q\in [1,+\infty]$ and $\Phi$ is a certain Young function.
\end{enumerate}

{In the rest of the paper, problem statement, basic assumptions, and construction of approximations for coercive or non-coercive Lyapunov functions are presented in Section~\ref{Problem setting}. The main results and their proofs of the $L^q$-ISS ($q\in [1,+\infty]$) in $L^1$-norm or weighted  $L^1$-norm for a class of nonlinear parabolic PDEs with Robin or Neumann or Dirichlet boundary disturbances are provided in Section~\ref{ISS L1}. The ISS in the setting of Orlicz spaces or Orlicz classes is addressed in Section~\ref{ISS Orlicz}. Some further comments and illustrative examples are presented in Section~\ref{Comments} and Section~\ref{Sec: examples}, respectively, followed by concluding remarks given in Section~\ref{Sec: Conclusion}. Well-posedness assessment for a class of quasi-linear parabolic equations and some technical development are presented in appendices.}

\emph{\textbf{Notations.}} In this paper, $\mathbb{R}_+ (\mathbb{R}_-)$ denotes the set of positive (negative) real numbers and $\mathbb{R}_{\geq 0} :=  \{0\}\cup\mathbb{R}_+$. Let $\Omega$ be a bounded domain in $\mathbb{R}^n (n\geq 1)$ with a sufficiently smooth boundary. We denote by $\partial \Omega$ and $\overline{\Omega}$ the boundary and the closure of $\Omega$, respectively. Denote by $|\Omega|$ (or $\partial \Omega$) the $n$-dimensional Lebesgue measure of $\Omega$ (or $(n-1)$-dimensional Lebesgue measure of $\partial\Omega$). Let $\mathbbm{d} =\max\limits_{x\in\partial\Omega}|x|$, which is the longest distance between the origin and $\partial\Omega$. Let $B_R(0)\subseteq \mathbb{R}^n$  be a ball with radius $R>0$ and center $0$. For any $T>0$, let $ Q_T=\Omega\times (0,T)$. 
Let $\mathcal {K}$, $\mathcal {L}$, and $\mathcal {K}\mathcal {L} $ be the standard classes of comparison functions as defined in, e.g., \cite{Mironchenko:2019b}.
Without special statements, the notations of function spaces and their norms used in this paper are standard which can be founded in, e.g., \cite{Evans:2010}.
\section{Problem setting and approximation of Lyapunov functions for ISS analysis}\label{Problem setting}
\subsection{Problem statement}
Given the following smooth functions:
\begin{align*}
&a,b_i,c,m_i,f,d\in C^{2}(\overline{\Omega}\times \mathbb{R}_{\geq 0}; \mathbb{R}),i=1,2,...,n,\\
&h\in C^{2}(\overline{\Omega}\times \mathbb{R}_{\geq 0}\times \mathbb{R}\times \mathbb{R}^n; \mathbb{R}),g\in C^{3}(\overline{\Omega}\times \mathbb{R}_{\geq 0}\times\mathbb{R}; \mathbb{R}),\\
&\psi \in C^{2}(\overline{\Omega}\times \mathbb{R}_{\geq 0}\times\mathbb{R}; \mathbb{R}),
w^0\in C^{2}(\overline{\Omega}; \mathbb{R}),
\end{align*}
and leting $\textbf{b}=(b_1,b_2,...,b_n), \textbf{m}=(m_1,m_2,...,m_n)$, we exploit the Lyapunov method to establish ISS and iISS estimates in different norms for the following nonlinear partial differential equation:
\begin{subequations}\label{LPE1'}
\begin{align}
& L_t[w](x,t)+N[w](x,t)
=f(x,t),  (x,t)\in\Omega\times \mathbb{R}_+,\\
&\mathscr{B}[w](x,t)
=d(x,t), (x,t)\in  \partial \Omega \times \mathbb{R}_+,\\
&w(x,0)=w^0(x), (x,t)\in  \Omega,
\end{align}
\end{subequations}
where $
 L_t[w](x,t)=w_t-\divi  \ (a\nabla w)
+\textbf{b}\cdot \nabla w +cw
$ and $N[w](x,t)=h(x,t,w,\nabla  w)+\bm{m}\cdot\nabla (g(x,t,w))$ are the linear and nonlinear parts of the equation, respectively, and
\begin{align}\label{Robin}
\mathscr{B}[w](x,t)=a\frac{\partial w}{\partial\bm{\nu}}+\psi(x,t,w),
\end{align}
or
\begin{align}\label{Dirichlet}
\mathscr{B}[w](x,t)=\psi(x,t,w),
\end{align}
represents the Neumann or nonlinear Robin, or nonlinear Dirichlet boundary condition, respectively, {with $\bm{\nu}$ denoting} the outer unit normal vector field on $\partial \Omega$. In general, $f$ and $d$ represent the distributed in-domain disturbance and boundary disturbance, respectively.

\begin{definition}\label{Def. ISS}
Let $p$ be a nonnegative weighting function. For $q,q_1,q_2\in [1,+\infty]$, System~\eqref{LPE1'} is said to be $L^{q}$-ISS in weighted $L^{q_1}$-norm w.r.t. boundary and in-domain disturbances $d$ and $f$ from, respectively, the $(L^{q}_{loc}(\mathbb{R}_+);L^{q_2}(\partial \Omega))$- and $(L^{q}_{loc}(\mathbb{R}_+);L^{q_2}(\Omega))$-spaces, if there exist functions $\beta_0\in \mathcal {K}\mathcal {L}$ and $ \gamma_0,\gamma_1\in \mathcal {K}$ such that the solution of \eqref{LPE1'} satisfies for any $T>0$:
\begin{align*}
    \|p^{\frac{1}{q_1}}w(\cdot,T)\|_{L^{q_1}(\Omega)}
    \leq& \beta_0\left(\|p^{\frac{1}{q_1}}w^0\|_{L^{q_1}(\Omega)},T\right)+\gamma_0 \left(\|d\|_{L^{q}((0,T);L^{q_2}(\partial\Omega))}\right)
      +\gamma_1\left(\|f\|_{L^{q}((0,T);L^{q_2}(\Omega))}\right).
\end{align*}
Furthermore, System~\eqref{LPE1'} is said to be exponential $L^{q}$-ISS in weighted $L^{q_1}$-norm w.r.t. boundary and in-domain disturbances $d$ and $f$ from, respectively, the $(L^{q}_{loc}(\mathbb{R}_+);L^{q_2}(\partial \Omega))$- and $(L^{q}_{loc}(\mathbb{R}_+);L^{q_2}(\Omega))$-spaces, if there are constants $M_0,\lambda>0$ such that $\beta_0( r,T)= M_0r\e^{-\lambda T}$ for all $T,r\geq 0$. Particularly, we omit the term ``weighted'' when $p\equiv 1$.
\end{definition}
\begin{definition} \label{Def. ISS'} Let $p$ be a nonnegative weighting function and $\Phi_1,\Phi_2$ be two $N$-functions\footnote{see \cite[Chapter 8]{AF2003} for a definition of ``$N$-function''.}. 
For $q\in [1,+\infty]$, System~\eqref{LPE1'} is said to be $L^q$-ISS in weighted $K_{\Phi_1}$-class w.r.t. boundary and in-domain disturbances $d$ and $f$ from, respectively, the $L^1_{loc}(\mathbb{R}_+$; $K_{\Phi_2}(\partial \Omega))$- and $L^1_{loc}(\mathbb{R}_+;K_{\Phi_2}(\Omega))$-classes, if there exist functions $\beta_0\in \mathcal {K}\mathcal {L}$ and $ \gamma_0,\gamma_1\in \mathcal {K}$ such that the solution of \eqref{LPE1'} satisfies for any $T>0$:
\begin{align*}
   \int_{\Omega}p\Phi_1(|w(x,T)|)\diff{x}
    \leq& \beta_0\bigg(\int_{\Omega}p\Phi_1(|w^0|)\diff{x},T\bigg)+\gamma_0(\mathcal {I}(d,\partial\Omega))+\gamma_1(\mathcal {I}(f,\Omega)),
   \end{align*}
   where for a certain set $\omega$ and a measurable function $v$, we denote $\mathcal {I}(v,\omega):=\left(\int_0^T\left(\int_{\omega}\Phi_2(|v|)\diff{x}\right)^q\diff{t}\right)^{\frac{1}{q}}$.
Furthermore, System~\eqref{LPE1'} is said to be exponential $L^q$-ISS in weighted $K_{\Phi_1}$-class w.r.t. disturbances $d$ and $f$ from, respectively, the $L^q_{loc}(\mathbb{R}_+;K_{\Phi_2}(\partial \Omega))$- and $L^q_{loc}(\mathbb{R}_+;K_{\Phi_2}(\Omega))$-classes, if there are constants $M_0,\lambda>0$ such that $\beta_0( r,T)= M_0r\e^{-\lambda T}$ for all $T,r\geq 0$. Particularly, we omit the term ``weighted'' when $p\equiv 1$.
\end{definition}
\begin{definition} \label{Def. ISS''}Let $\Phi_1,\Phi_2$ be two $N$-functions satisfying $\Delta_2$-condition\footnote{see \cite[Chapter 8]{AF2003} for a definition of         ``$\Delta_2$-condition''.}.
For $q\in [1,+\infty]$, System~\eqref{LPE1'} is said to be $L^q$-ISS in $L^{\Phi_1}$-norm w.r.t. boundary and in-domain disturbances $d$ and $f$ from, respectively, the $L^q_{loc}(\mathbb{R}_+$; $L^{\Phi_2}(\partial \Omega))$- and $L^q_{loc}(\mathbb{R}_+$; $L^{\Phi_2}(\Omega))$-spaces, if there exist functions $\beta_0\in \mathcal {K}\mathcal {L}$ and $ \gamma_0,\gamma_1\in \mathcal {K}$ such that the solution of \eqref{LPE1'} satisfies for any $T>0$:
\begin{align*}
    \|w(\cdot,T)\|_{L^{\Phi_1}(\Omega)}
    \leq& \beta_0\left(\|w^0\|_{L^{\Phi_1}(\Omega)},T\right)+\gamma_0 \left(\|d\|_{L^q((0,T);L^{\Phi_2}(\partial \Omega))}\right)
      +\gamma_1\left(\|f\|_{L^q((0,T);L^{\Phi_2}(\Omega))}\right).
\end{align*}
Furthermore, System~\eqref{LPE1'} is said to be exponential $L^q$-ISS in $L^{\Phi_1}$-norm w.r.t. boundary and in-domain disturbances $d$ and $f$ from, respectively, the $L^q_{loc}(\mathbb{R}_+;L^{\Phi_2}(\partial \Omega))$- and $L^q_{loc}(\mathbb{R}_+;L^{\Phi_2}(\Omega))$-spaces, if there are constants $M_0,\lambda>0$ such that $\beta_0( r,T)= M_0r\e^{-\lambda T}$ for all $T,r\geq 0$.
\end{definition}

\subsection{Basic assumptions}
We make the following structural assumptions.

\textbf{Assumption 1}\\
\textbf{(A1-1)}
 $\underline{a}\leq  a (x,t)\leq \overline{a}$ for all $(x,t)\in\overline{\Omega}\times \mathbb{R}_{\geq 0}
$,
where $\underline{a},\overline{a}$ are positive constants.\\
\textbf{(A1-2)} $\big(|h|+|h_t|+|h_s|+(1+|\eta|)\sum\limits_{i=1}^n|h_{\eta_i}|\big)\big|_{(x,t,s,\eta)}\leq \mu_1 (|s|)(1+|\eta|^2)
$ for all $(x,t,s,\eta)\in \overline{\Omega}\times \mathbb{R}_{\geq 0}\times \mathbb{R}\times \mathbb{R}^n$ with $\eta=(\eta_1,\eta_2,...,\eta_n)$, where $\mu_1\in C(\mathbb{R}_{\geq 0} ;\mathbb{R}_{\geq 0} ) $ is an increasing function.\\
\textbf{(A1-3)} $\big(\sum\limits_{i=1}^n|g_{x_i}|+|g_s|+|g_{ss}|\big)\big|_{(x,t,s)}\leq \mu_2(|s|)
$ for all $(x,t,s)\in\overline{ \Omega}\times \mathbb{R}_{\geq 0}\times \mathbb{R}$, where $\mu_2\in C(\mathbb{R}_{\geq 0} ;\mathbb{R}_{+} ) $ is an increasing function satisfying $ \int_{0}^{+\infty}\frac{1}{\mu_2(s)}\diff{s}=\infty$.\\
 \textbf{(A1-4)} $-\psi(x,t,s)s\leq \mu_3(x,t)(s^2+1)$ for all $(x,t,s)\in \partial{\Omega}\times \mathbb{R}_{\geq 0}\times \mathbb{R}$, where $\mu_3\in  C(\partial \Omega\times \mathbb{R}_{\geq 0};\mathbb{R})$.


\textbf{Assumption 2}\\
\textbf{(A2-1)} $
-h(x,t,s,\eta)s\leq \mu(x,t)s^2$ for all $(x,t,s,\eta)\in {\Omega}\times \mathbb{R}_{+}\times \mathbb{R}\times \mathbb{R}^n$, where $\mu\in C(\overline{\Omega}\times\mathbb{R}_{\geq 0} ;\mathbb{R}) $.\\
\textbf{(A2-2)}
 There exists a constant $s_0>0$ and a function $g_0\in C(\overline{ \Omega}\times \mathbb{R}_{\geq 0}; \mathbb{R})$ such that $|g(x,t,s)| \leq |g_0(x,t)||s| $ for all $(x,t,s)\in {\Omega}\times \mathbb{R}_{+}\times [-s_0,s_0]$.\\
  \textbf{(A2-3)}  There exists a constant $\underline{c}\geq 0$ such that
\begin{align}\label{c-}
\underline{c}\leq c-\divi \ \textbf{b}-\mu, \forall (x,t)\in{\Omega}\times \mathbb{R}_{+}.
\end{align}


\textbf{Assumption 3} For the problem with the Robin or Neumann boundary condition \eqref{Robin}:\\
\textbf{(A3-1)} $g(x,t,s)s\cdot\divi \ \bm{m}\leq 0$ for all $(x,t,s) \in  {\Omega}\times \mathbb{R}_{+} \times  \mathbb{R}$.\\
\textbf{(A3-2)} $\psi_0(x,t):=\inf\limits_{s\in\mathbb{R}\setminus\{0\}}\frac{\psi(x,t,s)
+s\bm{b}\cdot\bm{\nu}+g(x,t,s)\bm{m}\cdot\bm{\nu}}{s}$\
 \text{exists and}\ $\psi_0\geq \underline{\psi_0}$ for all $ (x,t)\in\partial\Omega\times\mathbb{R}_{+}$, where $\underline{\psi_0}\geq 0$ is a constant.


\textbf{Assumption 3'} For the problem with the Dirichlet boundary condition \eqref{Dirichlet}:\\
\textbf{(A3'-1)} $g$ and $\bm{m}$ satisfy $
g(x,t,s)s\cdot\divi \ (p\bm{m})\leq 0
$ for all $ (x,t,s)\in \overline{ \Omega}\times \mathbb{R}_{\geq 0}\times \mathbb{R}$, where $p\in C^2(\Omega)\cap C^1(\overline{\Omega})$ satisfies:
\begin{subequations}\label{condition p}
\begin{align}
 \divi \ (a\nabla p)+\bm{b}\cdot\nabla p&\leq -p_0, (x,t)\in  \Omega\times \mathbb{R}_{+},\label{condition pa}\\
p&=0,x\in  \partial \Omega,
\end{align}
\end{subequations}
with certain constant $p_0\geq 0$.\\
\textbf{(A3'-2)} $\psi(x,t,s)=\psi_1(x,t)\psi_2(s)$ for all $(x,t,s)\in \partial{\Omega}\times \mathbb{R}_{\geq 0}\times \mathbb{R}$, where $ \psi_1\geq \underline{\psi_1} >0$ for all $(x,t)\in \partial{\Omega}\times \mathbb{R}_{\geq 0}$ with some positive constant $ \underline{\psi_1}$, and $\psi_2$ is a strictly increasing and continuous function defined on $\mathbb{R}$.
%

We also impose the following compatibility conditions for the existence of a classical solution of system \eqref{LPE1'}.

\textbf{Assumption 4}\\
\textbf{(A4)} $
\mathscr{B}[w^0](x,t)=d(x,t) $ for all $ (x,t)\in \partial \Omega\times\{0\}$. 

\begin{remark} {}{Below are some remarks on the structural conditions.}

(i) System \eqref{LPE1'} contains the generalized Ginzburg-Landau equation if $h(x,t,w,\nabla w)\equiv c_1 w+c_2 w^3+c_3 w^5$ with $c_1\in \mathbb{R}$ and $c_2,c_3\in \mathbb{R}_{\geq 0}$  (see, e.g., \cite{Guo:1994}). Particularly, it contains the Chafee-Infante equation if $c_1\in \mathbb{R}_-,c_2\in \mathbb{R}_{\geq 0}$ and $c_3=0$.

(ii) System~\eqref{LPE1'} contains also the $n$-dimensional generalized Burgers equation if $\bm{m}\equiv (1,1,....,1)$ and $g(x,t,w)\equiv \frac{1}{q}w^q$ with a positive integer $q$ (see, e.g., \cite{Tersenov:2010,Zheng:201802}). For example, let $n=1$, $\Omega=(0,1),q=2$, $a\equiv 1,\bm{b}\equiv0,c\equiv1, h\equiv0$ and $\psi(x,t,s)= K(s+s^3)$ with $K\geq\frac{1}{2}$, System~\eqref{LPE1'} with the Robin boundary condition \eqref{Robin} becomes
    \begin{align*}
 &w_t+w_{xx}+w+ww_x=f(x,t),  (x,t)\in(0,1)\times \mathbb{R}_+,\\
&w_x(1,t)+K(w+w^3)
=d(1,t),t\in  \mathbb{R}_+,\\
& w_x(0,t)-K(w+w^3)
=-d(0,t), t\in  \mathbb{R}_+,\\
&w(x,0)=w^0(x), x\in  (0,1).
\end{align*}
Note that in \textbf{(A3-2)}:
\begin{align*}
\psi_0(i,t)&=\inf\limits_{s\in\mathbb{R}\setminus\{0\}}\frac{K(s+s^3)
+ \frac{1}{2}s^2\cdot(-1)^{i+1}}{s}\geq K-\frac{1}{2}\geq 0,i=0,1.
\end{align*}
Besides, it is easy to check that all the structural conditions in Assumption 2 and \textbf{(A3-1)} in Assumption 3 are satisfied.

(iii) Applying the comparison principle of elliptic equations to \eqref{condition p} (for fixed $t\in \mathbb{R}_+$), it follows that $p\geq 0$ in $\overline{\Omega}$. It is worthy noting that if $g(x,t,s)$ is odd in $s$ or $\bm{m}\equiv0$, there always exists $p\in C^2(\Omega)\cap C^1(\overline{\Omega})$ satisfying \eqref{condition p} so that \textbf{(A3'-1)} is satisfied trivially.

\end{remark}
Noting Assumption 1, 4 and \textbf{(A3'-2)}, one may verify that all the conditions of Theorem~\ref{well-posedness result} in Appendix~\ref{Appendix I} are satisfied for \eqref{LPE1'} with the boundary condition \eqref{Robin} or \eqref{Dirichlet}. Thus, we present first the following result on the well-posedness and omit the details of the proof.
\begin{proposition}\label{proposition existence}
Under the Assumption~1 and 4 (or Assumption~1, 4 and \textbf{(A3'-2)}), System~\eqref{LPE1'} with the Robin or Neumann boundary condition \eqref{Robin} (or  the Dirichlet boundary condition \eqref{Dirichlet}) admits a unique solution $w\in C^{2,1}(Q_T)\cap C(\overline{Q}_T)$  for any $T>0$.
\end{proposition}
\begin{remark} Below are some remarks on the assumptions for the well-posedness.

(i) {Assumption~1, 4 and \textbf{(A3'-2)} are only used for existence and regularity assessment of solutions to \eqref{LPE1'}. Particularly, Assumption~4 is used to assure that $w$ is continuous on the boundary of $ Q_T$, which is not necessary for ISS analysis. Indeed, Assumption~2, 3 (or 3') and \textbf{(A1-1)} are used for the establishment of ISS estimates in different norms or classes for System \eqref{LPE1'} with different boundary disturbances. }

(ii) Assumption~4 is imposed mainly for the following reasons:
\begin{enumerate}[$\bullet$]
  \item  {}{In general, in order to establish ISS estimates in different norms, one needs to define different solutions for PDE systems with different inputs. Since the aim of this paper is to introduce the {method of using approximations of Lyapunov functions for establishing} ISS estimates in different norms for nonlinear PDEs with inputs, e.g., $d,f$, and $w^0$, belonging to different spaces, for simplicity, we consider  only the classical solutions of \eqref{LPE1'} under Assumption 4, so that  every integration or derivative appearing in this paper makes a sense naturally.}

    \item {}{If $d,f$, and $w^0$ are smooth enough, one may remove Assumption~4 and consider a weak solution of \eqref{LPE1'} for {}{some} $T>0$ (with more {restrictions}, the weak solution can exist globally) in Sobolev spaces as shown in \cite{Amann:1986}. However, further arguments are needed for the existence of a weak solution in the setting of Orlicz spaces, which will be provided in our future work.}
\item {}{One may also remove Assumption 4 and consider a smooth solution  of \eqref{LPE1'} for the establishment of ISS estimates in different norms. In this case, more restrictions on the growth of $h,g$ in $s$ in \textbf{(A1-2)} and \textbf{(A1-3)} are needed as in \cite[Theorem 6.7, Chapter V]{Ladyzhenskaya:1968} for the existence of a smooth solution, which is continuous over $Q_T$, but not continuous on $\partial \Omega\times\{0\}$. Since no any other conditions are imposed on the growth of $h,g$ as in \cite[Theorem 6.7, Chapter V]{Ladyzhenskaya:1968}, the approach proposed in this paper {may be applied} to a wider class of nonlinear parabolic PDEs having a general form.}
\end{enumerate}

(iii) It should be mentioned that if $w^0$ and $d$ satisfy the condition in Assumption 4, the ``ISS estimates'' obtained in this paper are only in terms of such special $w^0$ and $d$ in a certain class as indicated in, e.g., \cite{Mironchenko:2019,karafyllis2017siam}, which are not the ISS properties in the sense that $w^0$ and $d$ should be arbitrary in certain spaces. Nevertheless, as mentioned above, it is still necessary to show how to apply the approximation of Lyapunov functions to establish ISS estimates in different norms, particularly in $L^1$-norm, for nonlinear PDEs by the Lyapunov method, which is the main goal of this paper. Besides, {we can consider the entropy solution or renormalized solution and establish ISS estimates  in $L^1$-norm for \eqref{LPE1'} with $w^0$ and $d$ from $L^1$-space as discussed in Section \ref{remarks existence}.}
\end{remark}

\subsection{Constructing approximations of Lyapunov functions}

To apply the Lyapunov method for ISS analysis, we define a convex function that will be used throughout this paper, i.e., for any $\tau\in \mathbb{R}_+$, let
\begin{align}\label{rho}
\rho_\tau(s)=\left\{\begin{aligned}
& |s|,\ |s|\geq \tau,\\
& -\frac{s^4}{8\tau^3}+\frac{3s^2}{4\tau}+\frac{3\tau}{8},\ |s|< \tau.
\end{aligned}\right.
\end{align}
It is easy to check that $\rho_\tau(s)$ is $C^2$-continuous in $s$ and satisfies for any $s\in \mathbb{R}$:
\begin{subequations}\label{property of rho}
\begin{align}
 \rho_\tau'(0)=0,0\leq |s|\leq\rho_\tau(s),|\rho_\tau'(s)|\leq 1,0\leq  \rho_\tau''(s),0\leq \rho_\tau(s)-\frac{3\tau}{8}\leq  \rho_\tau'(s)s\leq\rho_\tau(s) \leq |s|+\frac{3\tau}{8}.
\end{align}
\end{subequations}
{Fixing $\tau>0$, define $V_\tau(w):=\int_{\Omega}\rho_\tau(w)\diff{x}$. By \eqref{property of rho}, we find that
\begin{align*}
 \max\left\{\frac{3}{8}\tau|\Omega|,\|w\|_{L^1(\Omega)}\right\}\leq V_\tau(w)\leq \|w\|_{L^1(\Omega)} +\frac{3}{8}\tau|\Omega|.
\end{align*}
It is noted that $V_\tau(w)$ is a coercive functional since $V_\tau(w)\geq \|w\|_{L^1(\Omega)} $ and $V_\tau(w)\rightarrow +\infty $ as $ \|w\|_{L^1(\Omega)}\rightarrow +\infty$. However, $V_\tau(w)$ {cannot be taken as a Lyapunov function candidate} because $V_\tau(0) =\int_{\Omega}\rho_\tau(0)\diff{x} =\frac{3}{8}\tau|\Omega|>0$. On the other hand, it follows from the Lebesgue's Dominated Convergence Theorem that $V_\tau(w)\rightarrow  \int_{\Omega}|w|\diff{x}$ as $ \tau\rightarrow  0^+$ for any $w\in L^1(\Omega)$. Thus, $V_\tau(w)$ can be seen as an approximation of the coercive Lyapunov function candidate $\overline{V}(w):=\int_{\Omega}|w|\diff{x}$.


{Based on this idea, in this paper, we {construct} different approximations of Lyapunov functions {based on $\rho_\tau(\cdot)$ given in \eqref{rho}} for the establishment of ISS estimates in different norms or classes for the considered PDEs with different boundary disturbances from different spaces or classes. Specifically, letting $\tau>0$ be small enough and $\Phi$ be a certain Young function:
\begin{enumerate}[(i)]
\item for the $L^q$-ISS ($q\in [1,+\infty]$), or the $L^\Phi$-ISS, in $L^1$-norm of System \eqref{LPE1'} with Robin or Neumann boundary disturbances from $L^q$-space, or $L^\Phi$-space, we choose $ \int_{\Omega}\rho_\tau(w)\diff{x}$ as an approximation of the coercive Lyapunov function candidate $\int_{\Omega}|w|\diff{x}$;
\item for the $L^q$-ISS ($q\in [1,+\infty]$),  or the $L^\Phi$-ISS, in weighted $L^1$-norm (with a weighting function $p$) of System \eqref{LPE1'} with Dirichlet boundary disturbances from $L^q$-space, or $L^\Phi$-space, we choose $ \int_{\Omega}p\rho_\tau(w)\diff{x}$ as an approximation of the non-coercive Lyapunov function candidate $\int_{\Omega}p|w|\diff{x}$;
\item  for the $L^q$-ISS ($q\in [1,+\infty]$) in $L^\Phi$-norm (or $K_\Phi$-class) of System \eqref{LPE1'} with Robin or Neumann boundary disturbances from $L^\Phi$-space (or $K_\Phi$-class), we choose $ \int_{\Omega}\Phi(\rho_\tau(w))\diff{x}$ as an approximation of the coercive Lyapunov function candidate $\int_{\Omega}\Phi(|w|)\diff{x}$;
\item  for the $L^q$-ISS ($q\in [1,+\infty]$) in weighted $L^\Phi$-norm (or weighted $K_\Phi$-class, with a weighting function $p$) of System \eqref{LPE1'} with Dirichlet boundary disturbances from $L^\Phi$-space (or $K_\Phi$-class), we choose $ \int_{\Omega}p\Phi(\rho_\tau(w))\diff{x}$ as an approximation of the non-coercive Lyapunov function candidate $\int_{\Omega}p\Phi(|w|)\diff{x}$.
\end{enumerate}}


\section{$L^q$-ISS ($q\in [1,+\infty]$) in $L^1$-norm or weighted $L^1$-norm}\label{ISS L1}
\subsection{$L^q$-ISS in $L^1$-norm for nonlinear parabolic PDEs with Robin or Neumann boundary disturbances}
For System \eqref{LPE1'} with the Robin or Neumann boundary condition \eqref{Robin}, we have the following theorem.
\begin{theorem}\label{main result} Suppose that Assumption 1, 2, 3 and 4 hold. For any constant $q\in [1,+\infty]$ and $q'=\frac{q}{q-1}$, the following statements hold true:

(i) System \eqref{LPE1'} with the Robin or Neumann boundary condition \eqref{Robin} has the following $L^1$-estimate for all $T>0$:
 \begin{align}\label{042003}
\|w(\cdot,T)\|_{L^1( \Omega)}
\leq  &\e^{-\underline{c}T}\|w^0\|_{L^1( \Omega)}+
\| d\|_{L^1((0,T);L^1(\partial\Omega))}+\| f\|_{L^1((0,T);L^1(\Omega))}.
\end{align}
Furthermore, if $\underline{c}>0$ in \textbf{(A2-3)}, then System \eqref{LPE1'} with the Robin or Neumann boundary condition \eqref{Robin} is exponential $L^q$-ISS in $L^1$-norm w.r.t. boundary and in-domain disturbances $d$ and $f$ having the following estimate for all $T>0$:
\begin{align}\label{1803}
\!\!\|w(\cdot,T)\|_{L^1( \Omega)}
\leq  &\e^{-\underline{c}T}\|w^0\|_{L^1( \Omega)}
+ \bigg(\frac{1}{\underline{c}q'}\bigg)^{\frac{1}{q'}}
\| d\|_{L^q((0,T);L^1(\partial\Omega))}+ \bigg(\frac{1}{\underline{c}q'}\bigg)^{\frac{1}{q'}}\| f\|_{L^q((0,T);L^1(\Omega))}.
\end{align}

 (ii) Suppose that $\underline{c}=0$ in \textbf{(A2-3)} and $\underline{\psi_0}>0$ in $\bm{(\Psi_0})$. If
   there exist $j\in\{1,2,...,n\}$ and constants $ \underline{b}\in \mathbb{R},\overline{m}\in \mathbb{R}_{\geq 0}$ such that
\begin{align}
 &\underline{b}\leq \frac{b_j+a_{x_j}}{a},\forall (x,t)\in {\Omega}\times \mathbb{R}_{+},\label{AB}\\
& -m_jg(x,t,s)s\leq 0,\forall (x,t,s)\in {\Omega}\times \mathbb{R}_{+}\times\mathbb{R}\ (\text{or},
 \divi \ \bm{m}\neq 0\ \text{with}\ \bigg|\frac{m_j}{\divi \ \bm{m}}\bigg|\leq \overline{m},\forall (x,t)\in {\Omega}\times \mathbb{R}_{+}),\label{M}
   \end{align}
      then System \eqref{LPE1'} with the Robin or Neumann boundary condition \eqref{Robin} is exponential $L^q$-ISS in $L^1$-norm w.r.t. boundary and in-domain disturbances $d$ and $f$ having the following estimate for all $T>0$:
\begin{align}
\|w(\cdot,T)\|_{L^1( \Omega)}
\leq
 &C_{k,l}\e^{-\widehat{\underline{c}}T}\|{w^0}\|_{L^1( \Omega)}
+C_{k,l}\bigg(\frac{1}{\widehat{\underline{c}}q'}\bigg)^{\frac{1}{q'}}\!\!
\| d\|_{L^q((0,T);L^1(\partial\Omega))}+C_{k,l}\bigg(\frac{1}{\widehat{\underline{c}}q'}\bigg)^{\frac{1}{q'}}\| f\|_{L^q((0,T);L^1(\Omega))},
\label{1804}
\end{align}
  where
  \begin{subequations}\label{gains}
    \begin{align}
    \widehat{\underline{c}}=\frac{{}{\underline{a}}l\e^{-l\mathbbm{d}}}{k +\e^{l\mathbbm{d}}}\bigg(l+\underline{b}-\frac{2}{k}l\e^{l\mathbbm{d}}\bigg)>0,
    C_{k,l}=\frac{k+\e^{l\mathbbm{d} }}{k+\e^{-l\mathbbm{d} }}=1+\frac{\e^{-l\mathbbm{d}}+\e^{l\mathbbm{d}}}{k+ \e^{-l\mathbbm{d}}}>0,
    \end{align}
    \end{subequations}
     with any constants $l,k$ satisfying
    $l>\max\{0,-\underline{b}\}$ and $k>\max\{ \frac{2l\e^{l\mathbbm{d}}}{l+\underline{b}},\frac{\overline{a} l\e^{l\mathbbm{d}}}{\underline{\psi_0}}\}$ (or, $k>\max\{ \frac{2l\e^{l\mathbbm{d}}}{l+\underline{b}}$, $2l\e^{l\mathbbm{d}} \overline{m}$, $\frac{\overline{a} l\e^{l\mathbbm{d}}}{\underline{\psi_0}}\}$).
\end{theorem}

\begin{remark} \label{boundedness}
The parameters $ \widehat{\underline{c}}$ and $C_{k,l}$ in Theorem \ref{main result} (ii) are uniformly bounded in $l$ and $k$.
Indeed, letting $l\rightarrow+\infty$ (and thus $k\rightarrow+\infty$), we have $\widehat{\underline{c}}\rightarrow0$ and $C_{k,l}\rightarrow 1$. We deduce that there must be a positive constant $C_0$ suth that $0<\widehat{\underline{c}}<C_0$ and $1<C_{k,l}<1+C_0$ for all $l>\max\{0,-\underline{b}\}$ and $k>\max\{ \frac{2l\e^{l\mathbbm{d}}}{l+\underline{b}},\frac{\overline{a} l\e^{l\mathbbm{d}}}{\underline{\psi_0}}\}$ (or, $k>\max\{ \frac{2l\e^{l\mathbbm{d}}}{l+\underline{b}}$, $2l\e^{l\mathbbm{d}} \overline{m}$, $\frac{\overline{a} l\e^{l\mathbbm{d}}}{\underline{\psi_0}}\}$). Therefore, $C_{k,l}\cdot \big(\frac{1}{\widehat{\underline{c}}q'}\big)^{\frac{1}{q'}}$ is uniformly bounded in $l$ and $k$ with fixed $q\in[1,+\infty)$.
%

 %

\end{remark}
\begin{remark}
Compared with Theorem \ref{main result} (i), the positivity of $\underline{c}$ is weakened in Theorem~\ref{main result}~(ii). Therefore, in order to obtain an ISS estimate, it is natural to make more restrictions on $a$, $\bm{b},\bm{m}$ and $ \psi$. If $\psi_0(x,t)\equiv 0$, the result in Theorem~\ref{main result}~(ii) may fail to hold. For example, let $c=h=\psi =g=0,\bm{b}=\bm{m}=\bm{0}$ and $w^0$ be a non-zero constant, it is obvious that $w=w^0$ is the solution of \eqref{LPE1'} for $f=d=0$. However, $w=w^0$ does not tend to zero in any norm when $t\rightarrow +\infty$. Therefore, \eqref{LPE1'} is not ISS w.r.t. $d$ and $f$.
%
\end{remark}
\begin{pf*}{Proof of Theorem \ref{main result} (i):} 
For any $T>0$, let $w$ be the solution of the following equation:
\begin{subequations}\label{L1}
\begin{align}
 &L_t[w](x,t)+N[w](x,t)=f(x,t),(x,t)\in Q_T,\\
&\mathscr{B}[w](x,t)=d(x,t),(x,t)\in \partial \Omega\times (0,T),\\
& w(x,0)=w^0(x),x\in\Omega,
\end{align}
\end{subequations}
where $\mathscr{B}[w]$ is given in \eqref{Robin}.
{For any $\tau$ small enough (e.g., $0<\tau<s_0$), we choose $ \int_{\Omega}\rho_\tau(w)\diff{x}$ as an approximation of the coercive Lyapunov function candidate $\int_{\Omega}|w|\diff{x}$.} 
By direct computations, we have
\begin{align}\label{0301}
&\frac{\diff{}}{\diff{t}}\int_{\Omega}\rho_\tau(w)\diff{x}\notag\\
=&\!\int_{\partial \Omega}\! \!a\rho_\tau'(w)\frac{\partial w}{\partial\bm{\nu}}\diff{S}-\!
\int_{\Omega} \!\!a\rho_\tau''(w)|\nabla w|^2\diff{x}-\!\int_{\Omega}\!\!cw\rho_\tau'(w) \diff{x}-\int_{\Omega} \!\!(\textbf{b}  \cdot \nabla w)\rho_\tau'(w)\diff{x}-\int_{\Omega}\!\! (\bm{m}\cdot\nabla (g(x,t,w)))\rho_\tau'(w) \diff{x}\notag\\
&-\int_{\Omega}h(x,t,w,\nabla w)\rho_\tau'(w) \diff{x}+\int_{\Omega}f\rho_\tau'(w) \diff{x}.
\end{align}
We estimate the right hand side of \eqref{0301}. First, we have
\begin{align}\label{0302}
\int_{\partial \Omega} a\rho_\tau'(w)\frac{\partial w}{\partial\bm{\nu}}\diff{S}
=&\int_{\partial \Omega} \rho_\tau'(w)(d-\psi(x,t,w))\diff{S}
\leq \int_{\partial \Omega} |\rho_\tau'(w)||d|\diff{S}-\int_{\partial \Omega} \rho_\tau'(w)\psi(x,t,w)\diff{S}.
\end{align}
By \eqref{property of rho}, it follows that
\begin{align}\label{0903}
-\int_{\Omega}cw \rho_\tau'(w)\diff{x}
=& \int_{\Omega\cap\{-c\geq 0\}}-cw \rho_\tau'(w)\diff{x}+\int_{\Omega\cap\{-c< 0\}}-cw \rho_\tau'(w)\diff{x}\notag\\
\leq& \!\int_{\Omega\cap\{-c\geq 0\}}\!\!-c\rho_\tau(w)\diff{x}+\int_{\Omega\cap\{-c< 0\}}\!\!\!-c\bigg(\! \rho_\tau(w)-\frac{3}{8}\tau\!\bigg)\diff{x}\notag\\
=&-\int_{\Omega}c\rho_\tau(w)\diff{x}+\frac{3}{8}\tau\int_{\Omega\cap\{-c< 0\}}c\diff{x}\notag\\
\leq &-\int_{\Omega}c\rho_\tau(w)\diff{x}+\frac{3}{8}\tau\int_{\Omega} |c|\diff{x}.
\end{align}
We deduce from Green formulas and \eqref{property of rho} that
\begin{align}\label{1801}
-\int_{\Omega}(\textbf{b}  \cdot \nabla w)\rho_\tau'(w)\diff{x} =& -\int_{\Omega}\textbf{b}  \cdot \nabla (\rho_\tau(w))\diff{x}=\int_{\Omega}\rho_\tau(w)\divi \ \textbf{b}\diff{x}-\int_{\partial\Omega}\rho_\tau(w)\textbf{b}  \cdot \bm{\nu}\diff{S}\notag\\
\leq &\int_{\Omega}\rho_\tau(w)\divi \ \textbf{b}\diff{x}-\int_{\partial\Omega\cap\{ \textbf{b}  \cdot \bm{\nu}\geq 0\}}w\rho_\tau'(w)\textbf{b}  \cdot \bm{\nu}\diff{S}
-\int_{\partial\Omega\cap\{ \textbf{b}  \cdot \bm{\nu}< 0\}}\bigg(w\rho_\tau'(w)+\frac{3}{8}\tau\bigg)\textbf{b}  \cdot \bm{\nu}\diff{S}\notag\\
\!\!\!\leq &\int_{\Omega}\!\rho_\tau(w)\divi \ \textbf{b}\diff{x}-\!\int_{\partial\Omega}\!\!\!\!w\rho_\tau'(w)\textbf{b}  \cdot \bm{\nu}\diff{S}+\frac{3}{8}\tau\!\int_{\partial \Omega}\! \! |\textbf{b}|\diff{S}.
\end{align}
Similarly, it follows that
\begin{align}\label{2001}
-\int_{\Omega}h(x,t,w,\nabla w)\rho_\tau'(w) \diff{x}
= & \int_{\Omega\{w\neq 0\}}-h(x,t,w,\nabla w)w\cdot\frac{\rho_\tau'(w)}{w}  \diff{x}\notag\\
\leq &  \int_{\Omega\{w\neq 0\}}\mu w^2\cdot\frac{\rho_\tau'(w)}{w} \diff{x}\notag\\
\leq& \int_{\Omega}\mu\rho_\tau(w) \diff{x}+\frac{3}{8}\tau\int_{\Omega} |\mu|\diff{x}.
\end{align}
By integration by parts formula and noting \textbf{(A3-1)}, we have
\begin{align}\label{1301}
&-\int_{\Omega} (\bm{m}\cdot\nabla (g(x,t,w)))\rho_\tau'(w) \diff{x}\notag\\
=&-\int_{\partial \Omega} g(x,t,w)\rho_\tau'(w) \bm{m}  \cdot \bm{\nu}\diff{S}+\int_{\Omega} g(x,t,w)\rho_\tau''(w)\nabla w\cdot \bm{m}\diff{x}+\int_{\Omega} g(x,t,w)\rho_\tau'(w)\divi \ \bm{m}\diff{x}\notag\\
=& -\int_{\partial \Omega} g(x,t,w)\rho_\tau'(w) \bm{m}  \cdot \bm{\nu}\diff{S}
+\int_{\Omega} g(x,t,w)\rho_\tau''(w)\nabla w\cdot \bm{m}\diff{x}+\int_{\Omega\cap \{w\neq 0\}}\frac{1}{w^2}\cdot g(x,t,w)w\cdot\rho_\tau'(w)w\cdot\divi \ \bm{m}\diff{x}\notag\\
\leq&  -\int_{\partial \Omega} g(x,t,w)\rho_\tau'(w) \bm{m}  \cdot \bm{\nu}\diff{S}
+\int_{\Omega} |g(x,t,w)|\rho_\tau''(w)|\nabla w|| \bm{m}|\diff{x}\notag\\
\leq&
-\int_{\partial \Omega} g(x,t,w)\rho_\tau'(w) \bm{m}  \cdot \bm{\nu}\diff{S}
+\frac{1}{4\underline{a}}\int_{\Omega} |g(x,t,w)|^2\rho_\tau''(w)| \bm{m}|^2\diff{x}+\underline{a}\int_{\Omega} \rho_\tau''(w)| \nabla w|^2\diff{x},
\end{align}
with
\begin{align}\label{1302}
\int_{\Omega} |g(x,t,w)|^2\rho_\tau''(w)| \bm{m}|^2\diff{x}
\leq & \int_{\Omega} |g_{0}(x,t)|^2 |w|^2\rho_\tau''(w)| \bm{m}|^2\diff{x}\notag\\
\leq & \bigg(\max\limits_{ \overline{Q}_T} (|g_0|| \bm{m}|)^2\bigg)\cdot\int_{\Omega} |w|^2\rho_\tau''(w)\diff{x}\notag\\
\leq& \frac{3}{2\tau}\bigg(\max\limits_{ \overline{Q}_T} (|g_0| | \bm{m}|)^2\bigg)\cdot\int_{\Omega\cap\{|w|\leq \tau\}} |w|^2\bigg(1-\frac{w}{\tau}\bigg)\diff{x}\notag\\
\leq & \frac{3\tau}{2}|\Omega|\max\limits_{ \overline{Q}_T} (|g_0|| \bm{m}|)^2.
\end{align}
By  \textbf{(A3-2)} and \eqref{property of rho}, it follows that
 \begin{align}
 \int_{\partial \Omega}(\psi+w\bm{b}\cdot \bm{\nu}+g\bm{m}\cdot \bm{\nu})\rho_\tau'(w) \diff{S}
 =\int_{\partial\Omega\cap\{w\neq0\}}\rho_\tau'(w)w\cdot\frac{\psi+w\bm{b}\cdot \bm{\nu}+g\bm{m}\cdot \bm{\nu}}{w}\diff{S}
  \geq0.\label{1001}
 \end{align}
 It is obvious that
\begin{align}\label{0901}
\int_{ \Omega} \rho_\tau'(w)f\diff{x}\leq \int_{ \Omega} |f|\diff{x}.
\end{align}
From \eqref{0301} to \eqref{0901}, we get
\begin{align}\label{0902}
\frac{\diff{}}{\diff{t}}\int_{\Omega}\rho_\tau(w)\diff{x}
\leq& -\int_{\Omega}(c-\mu-\divi \ \textbf{b})\rho_\tau(w)\diff{x}
+\int_{\partial \Omega}|d|\diff{S}+\int_{\Omega}|f| \diff{x}\notag\\
&+\frac{3\tau}{8}\bigg(\int_{ \Omega}(|c|+|\mu|)\diff{x} +\frac{1}{\underline{a}}|\Omega|\max\limits_{ \overline{Q}_T} ((|g_0|| \bm{m}|)^2+| \bm{b}|^2)\bigg)\notag\\
 \leq& -\underline{c}\int_{\Omega}\rho_\tau(w)\diff{x}+A_1(t)+A_2(t)+B(t)\tau,
\end{align}
where
\begin{align}\label{Def. A(t)}
A_1(t)=\!\int_{\partial \Omega}|d|\diff{S},A_2(t)=\int_{\Omega}|f| \diff{x},
 B(t)=\!\frac{3}{8}\bigg(\!\!\int_{ \Omega}\!(|c|\!+\!|\mu|)\diff{x} +\!\!\int_{ \partial\Omega}\!\!\!|\bm{b}|\diff{S} \!+\!\frac{|\Omega|}{\underline{a}}\max\limits_{ \overline{Q}_T} (|g_0|| \bm{m}|)^2\!\bigg).
\end{align}
If $\underline{c}>0$, applying Gronwall's inequality to \eqref{0902} and using H\"{o}lder's inequality, we obtain
\begin{align}
&\int_{\Omega}\rho_\tau(w(x,T))\diff{x}\notag\\
\leq  &\e^{-\underline{c}T}\int_{\Omega}\rho_\tau(w^0)\diff{x}
+ \int_{0}^T( A_1(t)+A_2(t)+ B(t)\tau)\e^{-\underline{c}(T-t)}\diff{t}\label{A(t)}\\
\leq &\e^{-\underline{c}T}\int_{\Omega}\rho_\tau(w^0)\diff{x}
+ \left(\!\sum_{i=1}^2\!\|A_i\|_{L^q(0,T)}\!+\!\tau\|B\|_{L^q(0,T)}\!\right)\|\e^{-\underline{c}(T-\cdot)}\!\|_{L^{q'}(0,T)}\notag\\
\leq&\e^{-\underline{c}T}\int_{\Omega}\rho_\tau(w^0)\diff{x}
+ \left(\sum_{i=1}^2\|A_i\|_{L^q(0,T)}+\tau\|B\|_{L^q(0,T)}\right)\bigg(\frac{1}{\underline{c}q'}\bigg)^{\frac{1}{q'}},\label{041801}
\end{align}
where $q'=\frac{q}{q-1}$ with $q\in [1,+\infty]$. 
Letting $\tau\rightarrow 0$, we get \eqref{1803}. If $\underline{c}\geq 0$, then \eqref{042003} is a direct consequence of \eqref{A(t)}.
$\hfill\blacksquare$
\end{pf*}
\begin{pf*}{Proof of Theorem \ref{main result}~(ii):} 
Let $\beta(x)=k+\e^{l x_j}$, where $k, l$ are positive constants which will be chosen later.
Define the following quantities:
\begin{align*}
&\widehat{\bm{b}}= \bm{b}-\frac{2a}{\beta}\nabla \beta,\widehat{c}=c+\frac{1}{\beta}\bm{b}\cdot\nabla \beta-\frac{1}{\beta} \divi \ (a\nabla \beta),\notag\\
&
\widehat{h}(x,t,s,\eta)=\frac{1}{\beta}h(x,t,\beta s,\beta\eta +s\nabla \beta),\widehat{\mu}(x,t)=\mu(x,t),\\
&\widehat{\bm{m}}=\frac{1}{\beta}\bm{m},\widehat{g}(x,t,s)=g(x,t,\beta s),
\widehat{f}=\frac{1}{\beta}f, \widehat{d}=\frac{1}{\beta}d,\notag\\
&\widehat{w}^0=\frac{1}{\beta}w^0,
\widehat{\psi}(x,t,s)=\frac{1}{\beta} as(\nabla \beta \cdot \bm{\nu})+\frac{1}{\beta}\psi(x,t,\beta s),\\
&\widehat{I}(s)=\frac{\widehat{\psi}(x,t,s)+s\widehat{\bm{b}}\cdot \bm{\nu}+\widehat{g}(x,t,s)\bm{m}\cdot \bm{\nu}}{s},\notag\\
&
\widehat{L}_t[v]= v_t-\divi \ (a\nabla v)+\widehat{\bm{b}}\nabla v +\widehat{c}v,\notag\\
&\widehat{N}[v]=\widehat{h}(x,t,v,\nabla v)+\widehat{\bm{m}}\cdot\nabla (\widehat{g}(x,t,v)).
\end{align*}
 Putting $w=v\beta$ into \eqref{L1} and after some algebraic reduction, we obtain the following equation:
\begin{subequations}\label{hat}
\begin{align}
&\widehat{L}_t[v](x,t)+\widehat{N}[v](x,t) =\widehat{f}(x,t), (x,t)\in  Q_T,\\
&a\frac{\partial v}{\partial\bm{\nu}}+\widehat{\psi}(x,t,v)=\widehat{d}(x,t), (x,t)\in \partial \Omega\times (0,T),\\
&v(x,0)=\widehat{w}^0(x), x\in \Omega.
\end{align}
\end{subequations}
For the above system, it is easy to check that conditions \textbf{(A2-1)}, \textbf{(A2-2)} in Assumption 2 are satisfied (with a certain $s_0$ depending on $\beta$).

Based on the result of Theorem~\ref{main result}~(i), it suffices to find a constant $\widehat{\underline{c}} >0$ such that $\widehat{c}-\widehat{\mu}-\divi \ \widehat{\bm{b}}\geq\widehat{\underline{c}} >0$ satisfying \textbf{(A3-1)} and \textbf{(A3-2)} in Assumption~3 for System \eqref{hat}.

To this aim, we first choose $l>\max\{0,-\underline{b}\}$, $k> \frac{2l\e^{l\mathbbm{d}}}{l+\underline{b}}$ and deduce from \eqref{AB} that
\begin{align*}
\widehat{c}-\divi \ \widehat{\bm{b}}-\widehat{\mu}
=&c-\divi \ \bm{b}-{\mu}+\frac{1}{\beta}(\bm{b}\cdot\nabla \beta - \divi (a\nabla \beta))+2\divi \bigg(\frac{a\nabla \beta}{\beta}\bigg)\\
=&c-\divi \ \bm{b}-{\mu}+\frac{1}{\beta}(\bm{b}\cdot\nabla \beta + \divi (a\nabla \beta))-\frac{2a}{\beta^2}|\nabla \beta|^2\\
=& c-\divi \ \bm{b}-{\mu}+\frac{a}{\beta}l\e^{lx_j}\bigg(\frac{b_j+a_{x_j}}{a} +l-\frac{2}{\beta}l\e^{lx_j}\bigg)\notag\\
\geq &\frac{a}{\beta}l\e^{-l\mathbbm{d}}\bigg( l+\underline{b}-\frac{2}{\beta}l\e^{l\mathbbm{d}}\bigg)\notag\\
\geq &\frac{{\underline{a}}l\e^{-l\mathbbm{d}}}{k +\e^{l\mathbbm{d}}}\bigg(l+\underline{b}-\frac{2}{k}l\e^{l\mathbbm{d}}\bigg)\notag\\
:=& \widehat{\underline{c}}>0,\forall (x,t)\in{\Omega}\times \mathbb{R}_{+}.
\end{align*}
Note that according to \eqref{M},  if $ -m_jg(x,t,s)s\leq 0$ for all $(x,t,s)\in {\Omega}\times \mathbb{R}_{+}\times\mathbb{R}$, we deduce from \textbf{(A3-1)} and $\beta>0,l>0$ that
\begin{align}
 \widehat{g}(x,t,s)s\cdot \divi \ \widehat{\bm{m}}
 =g(x,t,\beta s)(\beta s) \cdot\frac{1}{\beta}\bigg( -\frac{l\e^{l x_j}}{\beta}m_j+ \divi \ {\bm{m}}\bigg)
 \leq 0,\forall(x,t,s) \in  { \Omega}\times\mathbb{R}_+\times \mathbb{R}.\label{1501}
\end{align}
If $\divi \ \bm{m}\neq 0$ with $|\frac{m_j}{\divi \ \bm{m}}|\leq \overline{m}$ for all $(x,t,s)\in {\Omega}\times \mathbb{R}_{+}$, we choose $\beta\geq 2l\e^{l\mathbbm{d}} \overline{m}$, which leads to $|-\frac{l\e^{lx_j}m_j}{\beta}\big|\leq \frac{1}{2}| \divi \ \bm{m}|$ for all $(x,t)\in {\Omega}\times \mathbb{R}_{+}$. It follows that 
 $\frac{1}{\beta}\big( -\frac{l\e^{l x_j}}{\beta}m_j+ \divi \ {\bm{m}}\big)$ and $\divi \ {\bm{m}}$ have the same sign, which, along with \textbf{(A3-1)}, guarantees that \eqref{1501} still hold.
Hence, the structural condition \textbf{(A3-1)} is satisfied for \eqref{hat}.

Now in order to verify \textbf{(A3-2)}, 
it suffices to prove that
$\widehat{I}(s)\geq 0 $ for all $(x,t,s)\in \partial{\Omega}\times \mathbb{R}_{+}\times \mathbb{R}\setminus\{0\}.$
Indeed, choosing $k>\frac{\overline{a} l\e^{l\mathbbm{d}}}{\underline{\psi_0}}$, it follows that
\begin{align*}
\widehat{I}(s)=&\frac{ \psi(x,t,\beta s) \! +\!(\beta s )\bm{b}\!\cdot \! \bm{\nu}\!+\!g(x,t,\beta s)  \bm{m}\!\cdot\! \bm{\nu}         }{\beta s}
 \!-\!\frac{a (\nabla \beta\!\cdot \!\bm{\nu})}{\beta}
\geq {\psi_0}-\frac{a l\e^{lx_{j}}}{k+\e^{l x_j}}
\geq\underline{\psi_0}-\frac{\overline{a} l\e^{l\mathbbm{d}}}{k}
\geq 0.
\end{align*}
By Theorem \ref{main result} (i), we have
\begin{align*}
\|v(\cdot,T)\|_{L^1( \Omega)}
\leq & \e^{-\widehat{\underline{c}}T}\|\widehat{w}^0\|_{L^1( \Omega)}
+\bigg(\frac{1}{\widehat{\underline{c}}q'}\bigg)^{\frac{1}{q'}}\cdot(
\| \widehat{d}\|_{L^q((0,T);L^1(\partial\Omega)}+\| \widehat{f}\|_{L^q((0,T);L^1(\Omega)}),
\end{align*}
for any constant $q\in [1,+\infty]$, where $q'=\frac{q}{q-1}$. Noting that $w=\beta v$, the desired result can be obtained.
$\hfill\blacksquare$
\end{pf*}

\subsection{$L^q$-ISS in weighted $L^1$-norm for nonlinear parabolic PDEs with Dirichlet boundary disturbances}
For System \eqref{LPE1'} with the Dirichlet boundary condition \eqref{Dirichlet}, we have the following theorem, which gives weighted $L^1$-estimates:
\begin{theorem}\label{main result 2}Suppose that Assumption 1, 2, 3' and 4 hold. For any constant $q\in [1,+\infty]$ and $q'=\frac{q}{q-1}$, the following statements hold true:

(i) {System \eqref{LPE1'} with the Dirichlet boundary condition \eqref{Dirichlet} has the following weighted $L^1$-estimate for all $T>0$:
\begin{align*}
\|pw(\cdot,T)\|_{L^1( \Omega)}
\leq & \e^{-\underline{c}T}\|pw^0\|_{L^1( \Omega)}+
\bigg \| a\psi^{-1}_2\left(\frac{d}{\psi_1}\right)\cdot|\nabla  p | \bigg\|_{_{L^1((0,T);L^1(\partial\Omega))}}+\| pf\|_{L^1((0,T);L^1(\Omega))},
\end{align*}
where $\psi_2^{-1}$ is the inverse of $\psi_2$ in \textbf{(A3'-2)}.}

{}{Furthermore, if $\underline{c}>0$ in \textbf{(A2-3)}, the following $L^1$-estimate holds for all $T>0$:
\begin{align}\label{ISS-Dirichlet}
\|pw(\cdot,T)\|_{L^1( \Omega)}
\leq  &\e^{-\underline{c}T}\|pw^0\|_{L^1( \Omega)}+\!\bigg(\!\frac{1}{\underline{c}q'}\!\bigg)^{\!\frac{1}{q'}}\!
\bigg \| a\psi^{-1}_2\!\!\left(\!\frac{d}{\psi_1}\!\right)\!\cdot\!|\nabla  p | \bigg\|_{_{L^q((0,T);L^1(\partial\Omega))}}+\bigg(\frac{1}{\underline{c}q'}\bigg)^{\frac{1}{q'}}\| pf\|_{L^q((0,T);L^1(\Omega))}.
\end{align}
Particularly, 
if $|\psi^{-1}_2(s)|\leq C|s|$ with a certain constant $C>0$ for all $s\in\mathbb{R}$ in \textbf{(A3'-2)}, then System \eqref{LPE1'} with \eqref{Dirichlet} is  $L^q$-ISS in weighted $L^1$-norm w.r.t. boundary and in-domain disturbances $d$ and $f$.}

(ii)
{}{Suppose that $\underline{c}=0$ in \textbf{(A2-3)} and $p_0>0$ in \textbf{(A3'-1)}. If there exist $j\in\{1,2,...,n\}$ and constants $ \underline{b}\in \mathbb{R},\overline{m}\in \mathbb{R}_{\geq 0}$  such that \eqref{AB} is satisfied and
\begin{align*}
&  -pm_jg(x,t,s)s\leq 0\ ,\forall (x,t,s)\in {\Omega}\times \mathbb{R}_{+}\times\mathbb{R}\ (\text{or},\divi \ (p\bm{m})\neq 0\ \text{with}\ \bigg|\frac{pm_j}{\divi \ (p\bm{m})}\bigg|\leq \overline{m},\forall (x,t)\in {\Omega}\times \mathbb{R}_{+}),
   \end{align*}
   then System \eqref{LPE1'} with the Dirichlet boundary condition \eqref{Dirichlet} has the following weighted $L^1$-estimate for all $T>0$:
\begin{align}\label{0413}
\|pw(\cdot,T)\|_{L^1( \Omega)}
\leq&  C_{k,l}\e^{-\widehat{\underline{c}}T}\|p{w^0}\|_{L^1(\Omega)}
+C_{k,l}\bigg(\frac{1}{\widehat{\underline{c}}q'}\bigg)^{\frac{1}{q'}}\bigg \| a\psi^{-1}_2\left(\frac{d}{\psi_1}\right)\cdot|\nabla  p | \bigg\|_{L^q((0,T);L^1(\partial\Omega))}\notag\\
&+C_{k,l}\bigg(\frac{1}{\widehat{\underline{c}}q'}\bigg)^{\frac{1}{q'}}\| pf\|_{L^q((0,T);L^1(\Omega))},
\end{align}
where $\widehat{\underline{c}}$ and $C_{k,l}$ are given by \eqref{gains} with any constants $l,k$ satisfying
    $l>\max\{0,-\underline{b}\}$ and $k>\max\{ \frac{2l\e^{l\mathbbm{d}}}{l+\underline{b}}$, $\frac{2\overline{a}l\e^{l\mathbbm{d} }}{p_0}\max\limits_{\overline{\Omega}}|p_j|\}$ (or $k>\max\{ \frac{2l\e^{l\mathbbm{d}}}{l+\underline{b}}$, $2l\e^{l\mathbbm{d}} \overline{m}$, $\frac{2\overline{a}l\e^{l\mathbbm{d}}}{p_0}\max\limits_{\overline{\Omega}}|p_j|\}$).
    Particularly, 
  if $|\psi^{-1}_2(s)|\leq C|s|$ with a certain constant $C>0$ for all $s\in\mathbb{R}$ in \textbf{(A3'-2)}, then System \eqref{LPE1'} with \eqref{Dirichlet} is  $L^q$-ISS in weighted $L^1$-norm w.r.t. boundary and in-domain disturbances $d$ and $f$.}
\end{theorem}
\begin{remark}
 As pointed out in Remark \ref{boundedness}, $C_{k,l}\cdot \big(\frac{1}{\widehat{\underline{c}}q'}\big)^{\frac{1}{q'}}$ in Theorem \ref{main result 2} (ii) is also uniformly bounded in $l$ and $k$ for fixed $q\in [1,+\infty)$.
 \end{remark}
\begin{pf*}{Proof of Theorem~\ref{main result 2}:} We consider \eqref{L1} with $\mathscr{B}[w]$ given by \eqref{Dirichlet}. 
As the proof can be proceeded in the same way as in that for Theorem~\ref{main result}, we present only the main steps.

For Theorem~\ref{main result 2}~(i), {we choose $ \int_{\Omega}p\rho_\tau(w)\diff{x}$ as an approximation of the non-coercive Lyapunov function candidate $\int_{\Omega}p|w|\diff{x}$.} Indeed, taking $p\rho_\tau'(w)$ as a test function and by direct computations, we have
\begin{align}\label{0301'}
\frac{\diff{}}{\diff{t}}\int_{\Omega}p\rho_\tau(w)\diff{x}
=&-\int_{ \partial\Omega} a\rho_\tau(w)\nabla  p \cdot \bm{\nu}\diff{S}+
\int_{\Omega} \rho_\tau(w)\divi \ (a\nabla p)\diff{x}
-
\int_{\Omega} ap\rho_\tau''(w)|\nabla w|^2\diff{x}+ \int_{\Omega} \rho_\tau(w)\divi \ ( p\bm{b})\diff{x}\
\notag\\
&-\int_{\Omega}cpw\rho_\tau'(w) \diff{x}+\int_{\Omega}pg(x,t,w) \rho_\tau''(w)\nabla w\cdot \bm{m}\diff{x}+\int_{\Omega}g(x,t,w) \rho_\tau'(w)\divi \ ( p\bm{m})\diff{x}\notag\\
&-\int_{\Omega}ph(x,t,w,\nabla w)\rho_\tau'(w) \diff{x}+\int_{\Omega}pf\rho_\tau'(w) \diff{x},
\end{align}
where 
\begin{align}
&\int_{\Omega} \rho_\tau(w)\divi \ ( p\bm{b})\diff{x}
=\int_{\Omega} \rho_\tau(w)\nabla p\cdot\bm{b}\diff{x}+\int_{\Omega} p\rho_\tau(w)  \divi \ \bm{b}\diff{x},  \label{0412b} \\
-&\!\int_{\Omega}\!\!cpw \rho_\tau'(w)\diff{x}
\leq \!-\!\!\int_{\Omega}\!\!cp\rho_\tau(w)\diff{x}\!+\!\frac{3}{8}\tau\!\!\int_{\Omega} \! \!|c|p\diff{x}, \label{0903'} \\
&\int_{\Omega}pg(x,t,w) \rho_\tau''(w)\nabla w\cdot \bm{m}\diff{x}\notag\\
\leq&
\int_{\Omega} p|g(x,t,w)|\rho_\tau''(w)| \bm{m}||\nabla w|\diff{x}\notag\\%
\leq & \frac{3\tau}{8\underline{a}}|\Omega|\max\limits_{ \overline{Q}_T} (p|g_0|^2| \bm{m}|^2)+\underline{a}\int_{\Omega}{p} \rho_\tau''(w)| \nabla w|^2\diff{x},\label{1302'} \\
&-\int_{\Omega}ph(x,t,w,\nabla w)\rho_\tau'(w) \diff{x}
\leq \int_{\Omega}\mu p\rho_\tau(w)  \diff{x}+\frac{3}{8}\tau\int_{\Omega} p|\mu|\diff{x},\label{2001'} \\
&\int_{ \Omega} p\rho_\tau'(w)f\diff{x}\leq \int_{ \Omega} p|f|\diff{x}.\label{0901'}
\end{align}
Note that by the boundary condition and \textbf{(A3'-2)}, we have $w=\psi^{-1}_2\big(\frac{d}{\psi_1}\big)$ on $\partial \Omega\times (0,T)$, which implies
\begin{align}\label{0412a}
-\int_{ \partial\Omega} a\rho_\tau(w)\nabla  p \cdot \bm{\nu}\diff{S}
\leq \int_{ \partial\Omega} a\cdot\bigg((\rho_\tau\circ\psi^{-1}_2)\left(\frac{d}{\psi_1}\right)\bigg)\cdot|\nabla  p |\diff{S}.
\end{align}
From \eqref{0301'} to \eqref{0412a}, and noting \textbf{(A3'-1)} and \textbf{(A2-3)}, we deduce that
\begin{align}\label{0301''}
\frac{\diff{}}{\diff{t}}\int_{\Omega}p\rho_\tau(w)\diff{x}
\leq&\underline{a}\int_{\Omega}p \rho_\tau''(w)| \nabla w|^2\diff{x}
-
\int_{\Omega} ap\rho_\tau''(w)|\nabla w|^2\diff{x}+\int_{\Omega} \rho_\tau(w)(\divi \ (a\nabla p)+\nabla  p\cdot\bm{b})\diff{x}
 \notag\\
&+\int_{\Omega}g(x,t,w) \rho_\tau'(w)\divi \ ( p\bm{m})\diff{x}-\int_{\Omega}p\rho_\tau(w)(c-\mu-\divi \ \bm{b})\diff{x}\notag\\
&+\int_{ \partial\Omega} a\cdot\bigg((\rho_\tau\circ\psi^{-1}_2)\left(\frac{d}{\psi_1}\right)\bigg)\cdot|\nabla  p |\diff{S}+\int_{\Omega}p|f| \diff{x}\notag\\
 &+\frac{3}{8}\tau\bigg(\int_{\Omega} p(|c|+|\mu|)\diff{x}+\frac{1}{\underline{a}}|\Omega|\max\limits_{ \overline{Q}_T} (p|g_0|^2| \bm{m}|^2)\bigg)\notag\\
\leq & -\underline{c}\int_{\Omega}p\rho_\tau(w)\diff{x}+A_1(t)+A_2(t)+B(t)\tau,
\end{align}
where
\begin{align*}
A_1(t)=&\int_{ \partial\Omega} a\cdot\bigg((\rho_\tau\circ\psi^{-1}_2)\left(\frac{d}{\psi_1}\right)\bigg)\cdot|\nabla  p |\diff{S},
A_2(t)=\int_{\Omega}p|f| \diff{x},\\
 B(t)=&\frac{3}{8}\bigg(\int_{\Omega} p(|c|+|\mu|)\diff{x}+\frac{1}{\underline{a}}|\Omega|\max\limits_{ \overline{Q}_T} (p|g_0|^2| \bm{m}|^2)\bigg).
\end{align*}
We prove Theorem \ref{main result 2} (i) only for the case of $\underline{c}>0$. Indeed, applying Gronwall's inequality to \eqref{0301''} and according to \eqref{041801}, we get \eqref{ISS-Dirichlet}. Moreover, it is obvious that the statement of Theorem~\ref{main result 2}~(i) holds true.

For Theorem \ref{main result 2} (ii), define $\beta$, $ \widehat{\bm{b}}$, $\widehat{c}$, $\widehat{\underline{c}}$, $\widehat{h}(x,t,s,\eta)$, $\widehat{\mu}$, $\widehat{\bm{m}}$, $\widehat{g}$,
$\widehat{f}$, $\widehat{w}^0$, $v$, $\widehat{L}_t[v]$, $\widehat{N}[v]
$ as in the proof of Theorem \ref{main result} (ii). Let $\widehat{d}=d$ and $\widehat{\psi}(x,t,s)=\psi_1(x,t)\widehat{\psi_2}(s)$ with $\widehat{\psi_2}( s)={\psi_2}(\beta s)$.
Then we get
\begin{subequations}\label{v-Dirichlet}
\begin{align}
&\widehat{L}_t[v](x,t)+\widehat{N}(x,t) =\widehat{f}, (x,t)\in  Q_T,\\
&\widehat{\psi}(x,t,v)=\widehat{d}, (x,t)\in \partial \Omega\times (0,T),\\
&v(x,0)=\widehat{w}^0(x), x\in \Omega.
\end{align}
\end{subequations}
Obviously, $ \psi_1\geq \underline{\psi_1} >0$ for all $(x,t)\in \partial{\Omega}\times \mathbb{R}_{\geq 0}$ with the constant $ \underline{\psi_1}$, and $\widehat{\psi_2}(s)$ is increasing and continuous in $s$.

Now choose $l>\max\{0,-\underline{b}\}$ and $k>\max\{ \frac{2l\e^{l\mathbbm{d}}}{l+\underline{b}},\frac{2\overline{a}l\e^{l\mathbbm{d} }}{p_0}\max\limits_{\overline{\Omega}}|p_j|\}$ (or $k>\max\{ \frac{2l\e^{l\mathbbm{d}}}{l+\underline{b}},2l\e^{l\mathbbm{d}} \overline{m},\frac{2\overline{a}l\e^{l\mathbbm{d} }}{p_0}\max\limits_{\overline{\Omega}}|p_j|\}$). Proceeding in the same way as in the proof of Theorem \ref{main result} (ii), we obtain
\begin{align*}
\widehat{g}(x,t,s)s\cdot\divi \ (p\widehat{\bm{m}})\leq 0,\forall (x,t,s) \in  {\Omega}\times \mathbb{R}_{+} \times  \mathbb{R},
\end{align*}
and
\begin{align*}
 \divi \ (a\nabla p)+\widehat{\bm{b}}\cdot\nabla p&\leq 0,   (x,t)\in  \Omega\times \mathbb{R}_{+},\\
p&=0,  x\in  \partial \Omega.
\end{align*}
Thus, all the conditions in Assumption~3' are satisfied. Applying Theorem~\ref{main result 2}~(i) to \eqref{v-Dirichlet}, we get
\begin{align*}
\|pv(\cdot,T)\|_{L^1( \Omega)}
\leq  &\e^{-\widehat{\underline{c}}T}\|pw^0\|_{L^1( \Omega)}+\bigg(\frac{1}{\widehat{\underline{c}}q'}\bigg)^{\frac{1}{q'}}\bigg \| a\widehat{\psi_2}^{-1}\left(\frac{\widehat{d}}{\psi_1}\right)\cdot|\nabla  p | \bigg\|_{L^q((0,T);L^1(\partial\Omega))}+\bigg(\frac{1}{\widehat{\underline{c}}q'}\bigg)^{\frac{1}{q'}}\| p\widehat{f}\|_{L^q((0,T);L^1(\Omega))},
\end{align*}
which implies \eqref{0413}. 
$\hfill\blacksquare$
\end{pf*}
\section{ISS and iISS in the setting of Orlicz spaces or Orlicz classes}\label{ISS Orlicz}
\subsection{Orlicz class $K_\Phi(\omega)$ and Orlicz space $L^\Phi(\omega)$}

In this section, we estimate solutions of \eqref{LPE1'} in the framework of Orlicz class or Orlicz space.
First of all, let $\phi\in C(\mathbb{R}_{\geq 0};\mathbb{R}_{\geq 0})\cap C^1(\mathbb{R}_{+};\mathbb{R}_{+})$ with $\phi(0)=0$ satisfy the Tolksdorf's condition:
\begin{align}\label{Tolksdorf}
 {\delta_0}\leq \frac{s\phi'(s)}{\phi(s)}\leq \delta_1,\ \forall\  s>0,
\end{align}
where $\delta_0,\delta_1>0$ are constants.

The structural condition of $\phi$ was firstly introduced by Tolksdorf in 1983 \cite{Tolksdorf:1983}, which is a natural generalization of the natural
conditions of Ladyzhenskaya and Ural'tseva in the existence and regularity theory of elliptic PDEs (see \cite{Lieberman:1991,Martinez:2008,Zheng:2019math,Zheng:2017math,ZhengGuo:2019,Zheng:2018Lyapunov}, etc.). Some typical examples of $\phi$ include:
(i) $\phi(s)=s^{q-1},\forall q>1$;
(ii) $\phi(s)=\ln(1+c_1s)+c_2s,\ \ \forall c_1,c_2>0$;
 (iii) $\phi(s)=(\ln (s+c_1))^{c_2}s^{q-1} ,\forall c_1\geq \e,c_2>0,q> 1$.
More examples can be found in \cite{ZhengGuo:2019}.

Let $\Phi(s)=\int_{0}^s\phi(\tau)\diff{\tau}$. Under the assumption \eqref{Tolksdorf}, $\Phi$ has the following properties.
\begin{lemma} 
\label{property of phi} The following results hold true:
\begin{enumerate}[(i)]
\item $\min\{k^{\delta_0},k^{\delta_1}\}\phi(s)\leq \phi(ks)\leq \max\{k^{\delta_0},k^{\delta_1}\}\phi(s)$, $\forall  k,s\geq
 0$;
\item $\Phi$ is $C^{2}$-continuous on $(0,+\infty)$ and convex on $[0,+\infty)$;
\item $\frac{1}{1+\delta_1}s\phi(s)\leq \Phi(s)\leq \frac{1}{1+\delta_0}s\phi(s)$,  $\forall  s\geq
 0$;
 \item $\frac{1+\delta_0}{1+\delta_1} \min\{k^{1+\delta_0},k^{1+\delta_1}\}\Phi(s)$ $\leq$ $\Phi(ks)$ $\leq$ $\frac{1+\delta_1}{1+\delta_0}\times$ $\max\{k^{1+\delta_0},k^{1+\delta_1}\}\Phi(s)$, $\forall  k,s\geq 0$.
\end{enumerate}
\end{lemma}
\begin{pf*}{Proof:} It is obvious that (ii) holds true. The proofs of (i) and (iii) are provided in \cite{Lieberman:1991} and \cite{ZhengGuo:2019}, respectively. Finally, (iv) is a consequence of (i) and (iii).
$\hfill\blacksquare$
\end{pf*}

Obviously, $\Phi$ is a $N$-function and satisfies {$ \Delta_2$-condition}. Associated with the $N$-function $\Phi$, we recall the definitions of Orlicz class and Orlicz space, which can be found in \cite[Chapter 8]{AF2003}. Let $ \omega$ be an bounded domain in $\mathbb{R}^{n} $ or $\mathbb{R}^{n-1}$ or $\mathbb{R}$. The Orlicz class $ K_\Phi(\omega)$ is the set of all measurable functions $u$ defined on $\omega$ that satisfy $\int_\omega\Phi(|u(y)|)\diff{y}<+\infty$. The Orlicz space $L^{\Phi}(\omega)$ is the linear hull of the Orlicz class $K_\Phi(\omega) $, which is a Banach space w.r.t. the Luxemburg  norm defined on $L^{\Phi}(\omega)$, i.e., $\|u\|_{L^{\Phi}(\omega)}=\inf\left\{k>0;\ \int_{\omega}\Phi\left(\frac{|u(y)|}{k}\right)\diff{y}\leq
1\right\}$. 
It is worthy noting that $L^{\Phi}(\Omega)=L^{q}(\Omega)$ for $\delta_1=\delta_0=q-1,q\in[1,+\infty)$. As $\phi$ is strictly increasing, we can define its inverse function $\phi^{-1}$. Let $\widetilde{\Phi}(s)=\int^s_0\phi^{-1}(\tau)\diff{\tau}$ for $s\geq 0$, which is the complementary $N$-function of $\Phi(s)$.
\subsection{$L^\Phi$-ISS in $L^1$-norm or weighted $L^1$-norm for nonlinear parabolic PDEs with boundary disturbances}
\begin{theorem}\label{main result Phi}
For System \eqref{LPE1'} with the Robin or Neumann boundary condition \eqref{Robin}, the following statements hold true:

(i) Under the same assumptions as in Theorem \ref{main result} (i), if $\underline{c}>0$ in \textbf{(A2-3)}, then System \eqref{LPE1'} with the Robin or Neumann boundary condition \eqref{Robin} is exponential $L^{\Phi}$-ISS in $L^1$-norm w.r.t. boundary and in-domain disturbances $d$ and $f$ having the following estimate for all $T>0$:
\begin{align}
\|w(\cdot,T)\|_{L^1( \Omega)}
\leq  &\e^{-\underline{c}T}\|w^0\|_{L^1( \Omega)}+ C(\underline{c},\delta_0,\delta_1 )
\| d\|_{L^\Phi((0,T);L^1(\partial\Omega))}+ C(\underline{c},\delta_0,\delta_1 )\| f\|_{L^\Phi((0,T);L^1(\Omega))},\label{4.2}
\end{align}
where $C(\underline{c},\delta_0,\delta_1 )=2\max\limits_{i\in\{0,1\}}\left\{\left(\frac{1}{\underline{c}}\frac{\delta_1^2}{\delta_0(1+\delta_1)} \widetilde{\Phi}(1)\right)^{\frac{1}{1+\delta_i}}\right\}$ is a positive constant.

(ii) Under the same assumptions as in Theorem \ref{main result} (ii), System \eqref{LPE1'} with the Robin or Neumann boundary condition \eqref{Robin} is exponential $L^{\Phi}$-ISS in $L^1$-norm w.r.t. boundary and in-domain disturbances $d$ and $f$ having the following estimate for all $T>0$:
\begin{align}
\|w(\cdot,T)\|_{L^1( \Omega)}
\leq  C_{k,l}\e^{-\widehat{\underline{c}}T}\|{w^0}\|_{L^1( \Omega)}+
C_{k,l}C'(\widehat{\underline{c}},\delta_0,\delta_1 )
\| d\|_{L^\Phi((0,T);L^1(\partial\Omega))}+C_{k,l}C'(\widehat{\underline{c}},\delta_0,\delta_1 )\| f\|_{L^\Phi((0,T);L^1(\Omega))},\label{4.3}
\end{align}
where $\widehat{\underline{c}},C_{k,l}$ are the same as in Theorem \ref{main result} (ii), $C'(\widehat{\underline{c}},\delta_0,\delta_1 )=2\max\limits_{i\in\{0,1\}}\bigg\{\left(\frac{1}{\widehat{\underline{c}}}\frac{\delta_1^2}{\delta_0(1+\delta_1)} \widetilde{\Phi}(1)\right)^{\frac{1}{1+\delta_i}}\bigg\}$ is a positive constant.
\end{theorem}

\begin{theorem}\label{main result 2 Phi}
For System \eqref{LPE1'} with the Dirichlet boundary condition \eqref{Dirichlet}, the following statements hold true:

 (i) {}{Under the same assumptions as in Theorem \ref{main result 2} (i), if $\underline{c}>0$ in \textbf{(A2-3)}, then System \eqref{LPE1'} with the Dirichlet boundary condition \eqref{Dirichlet} has the following weighted $L^1$-estimate for all $T>0$:
\begin{align*}
\|pw(\cdot,T)\|_{L^1( \Omega)}
\leq  \e^{-\underline{c}T}\|pw^0\|_{L^1( \Omega)}
+C(\underline{c},\delta_0,\delta_1 )
\bigg \| a\psi^{-1}_2\left(\frac{d}{\psi_1}\right)\cdot|\nabla  p | \bigg\|_{_{L^\Phi((0,T);L^1(\partial\Omega))}}
+C(\underline{c},\delta_0,\delta_1 )\| pf\|_{L^\Phi((0,T);L^1(\Omega))},
\end{align*}
where $\psi_2^{-1}$ is the inverse of $\psi_2$ in \textbf{(A3'-2)}, and $C(\underline{c},\delta_0,\delta_1 )$ is the same as in Theorem \ref{main result Phi} (i).} {Furthermore, if $|\psi^{-1}_2(s)|\leq C|s|$ with a certain constant $C>0$ for all $s\in\mathbb{R}$ in \textbf{(A3'-2)}, then System \eqref{LPE1'} with \eqref{Dirichlet} is  $L^\Phi$-ISS in $L^1$-norm w.r.t. boundary and in-domain disturbances $d$ and $f$.}

(ii) {}{Under the same assumptions as in Theorem \ref{main result 2} (ii), System \eqref{LPE1'} with the Dirichlet boundary condition \eqref{Dirichlet} has the following weighted $L^1$-estimate for all $T>0$:
\begin{align*}
\|pw(\cdot,T)\|_{L^1( \Omega)}
\leq & C_{k,l}\e^{-\widehat{\underline{c}}T}\|p{w^0}\|_{L^1(\Omega)}
+C_{k,l}C'(\widehat{\underline{c}},\delta_0,\delta_1 )
\bigg \| a\psi^{-1}_2\left(\frac{d}{\psi_1}\right)\cdot|\nabla  p | \bigg\|_{_{L^\Phi((0,T);L^1(\partial\Omega))}}\notag\\
&
+C_{k,l}C'(\widehat{\underline{c}},\delta_0,\delta_1 )\| pf\|_{L^\Phi((0,T);L^1(\Omega))},
\end{align*}
where $\widehat{\underline{c}},C_{k,l}$ are the same as in Theorem~\ref{main result 2}~(ii), and $C'(\widehat{\underline{c}},\delta_0,\delta_1 )$ is the same as in Theorem~\ref{main result Phi}~(ii).}
\end{theorem}

We only prove Theorem~\ref{main result Phi}~(i). The proof of Theorem~\ref{main result Phi}~(ii) and Theorem~\ref{main result 2 Phi} can be proceeded in a similar way.
\begin{pf*}{Proof of Theorem~\ref{main result Phi}~(i):} We choose $ \int_{\Omega}\rho_\tau(w)\diff{x}$ as an approximation of the coercive Lyapunov function candidate $\int_{\Omega}|w|\diff{x}$. It suffices to deduce from \eqref{A(t)} and H\"{o}lder's inequality (Lemma \ref{property of inverse of phi'} (iii)) that
\begin{align}
\int_{\Omega}\rho_\tau(w(x,T))\diff{x}
\leq  &\e^{-\underline{c}T}\!\!\int_{\Omega}\!\rho_\tau(w^0)\diff{x}
+\!\!\int_{0}^T\!\!\!\!\left(\!\sum\limits_{i=1}^2 \!A_i(t)+ B(t)\tau\!\!\right)\!\e^{-\underline{c}(T-t)}\!\diff{t}\notag\\
\leq &2\!\left(\!\sum\limits_{i=1}^2\|A_i\|_{L^{\Phi}(0,T)}+\tau\|B\|_{L^{\Phi}(0,T)}\!\right)\!\|\e^{-\underline{c}(T-\cdot)}\|_{L^{\widetilde{\Phi}}(0,T)}
+ \e^{-\underline{c}T}\int_{\Omega}\rho_\tau(w^0)\diff{x}
,\label{042002}
\end{align}
where $A_1(t),A_2(t)$ and $B(t)$ are the same as  in \eqref{Def. A(t)}.

For $t\in [0,T]$, it follows from Lemma \ref{property of inverse of phi} (i) that
\begin{align*}
\widetilde{\Phi}(\e^{-\underline{c}(T-t)})\leq  \frac{\delta_1(1+\delta_0)}{\delta_0(1+\delta_1)} \e^{-\underline{c}(T-t)\big(1+\frac{1}{\delta_1}\big)}\widetilde{\Phi}(1),
\end{align*}
which gives
\begin{align*}
\int_0^T\!\!\widetilde{\Phi}(\e^{-\underline{c}(T-t)})\diff{t}\leq  & \frac{\delta_1(1+\delta_0)}{\delta_0(1+\delta_1)} \widetilde{\Phi}(1)\!\int_0^T\!\!\e^{-\underline{c}(T-t)\big(1+\frac{1}{\delta_1}\big)}\!\diff{t}
\leq \frac{\delta_1^2(1+\delta_0)}{\underline{c}\delta_0(1+\delta_1)^2} \widetilde{\Phi}(1),
\end{align*}

 Noting  Lemma \ref{property of inverse of phi'} (ii) and applying Lemma \ref{property of inverse of phi'} (i) to $\widetilde{\Phi}$, we claim that
  \begin{align}
\|\e^{-\underline{c}(T-\cdot)}\|_{L^{\widetilde{\Phi}}(0,T)}\leq \max\limits_{i\in\{0,1\}}\bigg\{\bigg(\frac{1}{\underline{c}}\frac{\delta_1^2\widetilde{\Phi}(1)}{\delta_0(1+\delta_1)} \bigg)^{\frac{1}{1+\delta_i}}\bigg\}.\label{042001}
\end{align}

By \eqref{042002}, \eqref{042001} and letting $\tau\rightarrow 0$, we get \eqref{4.2}.
$\hfill\blacksquare$
\end{pf*}

\subsection{$L^q$-ISS ($q\in [1,+\infty]$) in $L^\Phi$-norm or weighted $K_\Phi$-class  for nonlinear parabolic PDEs with boundary disturbances}
In this subsection, for simplicity, we only consider ISS and $L^1$-ISS of System \eqref{LPE1'} for $g\equiv 0$ and inputs $w^0,d,f$ from $K_\Phi$-class or $L^\Phi$-space. Without special statements, we replace 
\textbf{(A2-3)} in Assumption~2 and \textbf{(A3-2)} in Assumption~3 by, respectively, the following conditions:\\
\textbf{(A2-3')} There exists a constant $\underline{c}> 0$ such that
\begin{align}\label{c-''}
\!\!\!\underline{c}\leq c-\mu-\max\bigg\{\frac{\divi \ \textbf{b} }{1+\delta_0},\frac{\divi \ \textbf{b}}{1+\delta_1}\bigg\},
 \forall (x,t)\in{\Omega}\times \mathbb{R}_{+};
\end{align}
\textbf{(A3-2')}
$\psi_0(x,t):=\inf\limits_{s\in\mathbb{R}\setminus\{0\}}\frac{\psi(x,t,s)
+s\min\big\{\frac{\bm{b}\cdot \nu }{1+\delta_0},\frac{\bm{b}\cdot \nu }{1+\delta_1}\big\}}{s}$\
 \text{exists and}\ $\psi_0\geq \underline{\psi_0}$ for all $ (x,t)\in\partial\Omega\times\mathbb{R}_{+}$, where $\underline{\psi_0}> 0$ is a constant.

For $\varepsilon>0$ small enough, let $\lambda:=\underline{c}(1+\delta_0)-\varepsilon>0$ and $C(\varepsilon):=\frac{\delta_1(1+\delta_0)}{\delta_0(1+\delta_1)} \bigg(\frac{1+\delta_0}{1+\delta_1} \varepsilon\bigg)^{-\frac{1+\delta_0}{\delta_0(1+\delta_1)}}$.

For System \eqref{LPE1'} with the Robin or Neumann boundary condition \eqref{Robin}, we have the following result.
\begin{theorem}\label{main result'}
{}{Suppose that Assumption 1, 2, 3, and 4 hold. For any $q\in [1,+\infty]$, System \eqref{LPE1'} (for $g\equiv0$) with the Robin or Neumann boundary condition \eqref{Robin} is exponential $L^q$-ISS in $K_\Phi$-class w.r.t. boundary and in-domain disturbances $d$ and $f$ from, respectively, the $L^q_{loc}(\mathbb{R}_+;K_{\Phi}(\partial \Omega))$- and $L^q_{loc}(\mathbb{R}_+;K_{\Phi}(\Omega))$-classes respectively, having the following estimate for all $T>0$:
 \begin{align}
\int_{\Omega}\!\Phi(|w(x,T)|)\diff{x}
\leq & \e^{-\lambda T}\!\!\!\int_{\Omega}\!\!\Phi(|w^0|)\diff{x}\!+\!C(\varepsilon)\bigg(\! \frac{1}{\lambda q'}\!\bigg)^{\!\!\frac{1}{q'}}\!\!\bigg(\!\!\int_0^T\!\!\!\!\bigg(\!\int_{\partial\Omega}\!\!\!\Phi(|d|)  \!\diff{S}\!\bigg)^{\!\!q}\!\!\diff{t}\!\bigg)^{\!\!\frac{1}{q}}\notag\\
&+C(\varepsilon)\bigg( \!\frac{1}{\lambda q'}\!\bigg)^{\!\frac{1}{q'}}\!\!\bigg(\!\int_0^T\!\!\!\bigg(\!\int_{\Omega}\!\!\Phi(|f|)  \diff{S}\bigg)^{q}\!\!\diff{t}\bigg)^{\!\frac{1}{q}},\label{041802'}
\end{align}
which, along with Lemma~\ref{property of inverse of phi'}~(i), implies further that System~\eqref{LPE1'} with the Robin or Neumann boundary condition \eqref{Robin} is exponential $L^q$-ISS in $L^\Phi$-norm w.r.t. boundary and in-domain disturbances $d$ and $f$ from  $L^q_{loc}(\mathbb{R}_+;K_{\Phi}(\partial \Omega))$- and $L^q_{loc}(\mathbb{R}_+;K_{\Phi}(\Omega))$-classes respectively.}
\end{theorem}
%

 For System \eqref{LPE1'} with the Dirichlet boundary condition \eqref{Dirichlet}, we have the following result.
\begin{theorem}\label{main result 2'}
{}{Suppose that Assumption 1, 2, 3' and 4 hold. For any $q\in [1,+\infty]$, System \eqref{LPE1'} (for $g\equiv 0$) with the Dirichlet boundary condition \eqref{Dirichlet} has the following weighted $K_\Phi$-estimate for all $T>0$:
\begin{align}\label{ISS-Dirichlet''}
\int_{\Omega}p\Phi(|w(x,T)|)\diff{x}
\leq &\bigg(\! \frac{1}{\lambda q'}\!\bigg)^{\!\!\frac{1}{q'}} \!\!\bigg\{\!\!\bigg(\!\int_0^T\!\!\!\bigg(\!\int_{ \partial\Omega} \! \!\!a\!\cdot\!\bigg(\!(\Phi\circ\psi^{-1}_2)\!\left(\frac{|d|}{\psi_1}\!\right)\!\!\bigg)\!\cdot\!|\nabla  p |\diff{S}\bigg)^{\!\!q}\!\!\diff{t}\!\bigg)^{\!\!\frac{1}{q}}\notag\\
&+C(\varepsilon)\frac{1+\delta_1}{1+\delta_0} \! \bigg(\!\int_0^T\!\!\!\!\bigg(\!\!\int_{\Omega}\!\max\limits_{i\in\{0,1\}}\{ |p|^{1+\delta_i}\}\Phi(|f|)  \diff{x}\!\bigg)^{\!\!q}\!\!\!\diff{t}\!\bigg)^{\!\!\frac{1}{q}}\bigg\}+\e^{-\lambda T}\int_{\Omega} p\Phi(|w^0|)\diff{x}.
\end{align}
Furthermore, if $|\psi^{-1}_2(s) |\leq C|s|$ on $\mathbb{R}$ with a certain constant $C>0$, then
System \eqref{LPE1'} with \eqref{Dirichlet} is  exponential $L^q$-ISS in weighted $K_{\Phi}$-class w.r.t. boundary and in-domain disturbances $d$ and $f$ from $L^q_{loc}(\mathbb{R}_+;K_{\Phi}(\partial \Omega))$- and $L^q_{loc}(\mathbb{R}_+;K_{\Phi}(\Omega))$-classes respectively.}
\end{theorem}
\begin{pf*}{Proof of Theorem~\ref{main result'}:} We consider \eqref{L1} with $\mathscr{B}[w]$ given by \eqref{Dirichlet}. We proceed in the same way as in the proof of Theorem~\ref{main result} and only present the main steps.

{We choose $ \int_{\Omega}\Phi(\rho_\tau(w))\diff{x}$ as an approximation of the coercive Lyapunov function candidate $\int_{\Omega}\Phi(|w|)\diff{x}$.} Indeed,
taking $\phi(\rho_\tau(w))\rho_\tau'(w)$ as a test function and by direct computations, we have
\begin{align}\label{1500'}
\frac{\diff{}}{\diff{t}}\int_{\Omega}\Phi(\rho_\tau(w))\diff{x}
=&\int_{ \partial\Omega}\! a\phi(\rho_\tau(w))\rho_\tau'(w)\nabla  w \!\cdot\! \bm{\nu}\diff{S}\!-\!
\int_{\Omega} \!a \bigg(\phi'(\rho_\tau(w))(\rho_\tau'(w))^2+\phi(\rho_\tau(w))\rho_\tau''(w)\bigg)|\nabla w|^2\diff{x}\notag\\
&-\int_{\Omega}c\phi(\rho_\tau(w))w\rho_\tau'(w)\diff{x}
{-\int_{\partial\Omega}\Phi(\rho_\tau(w)) \bm{b}\cdot \bm{\nu}\diff{S}+\int_{\Omega}\Phi(\rho_\tau(w))\divi \ \bm{b}\diff{x}}
\notag\\
&-\int_{\Omega}h(x,t,w,\nabla w)\phi(\rho_\tau(w))\rho_\tau'(w) \diff{x}+\int_{\Omega}f\phi(\rho_\tau(w))\rho_\tau'(w) \diff{x}.
\end{align}
By \eqref{property of rho}, Lemma \ref{property of phi} (iii) and Lemma \ref{property of inverse of phi} (ii, iii), it follows that
\begin{align}\label{1506}
\int_{ \partial\Omega} a\phi(\rho_\tau(w))\rho_\tau'(w)\nabla  w \cdot \bm{\nu}\diff{S}
=&\int_{ \partial\Omega}\!\! \phi(\rho_\tau(w))\rho_\tau'(w)d\diff{S}\!-\!\!\int_{ \partial\Omega} \!\!\phi(\rho_\tau(w))\rho_\tau'(w)\psi(x,t,w)\diff{S}\notag\\
\leq&\int_{ \partial\Omega} \phi(\rho_\tau(w))|d|\diff{S}-\int_{ \partial\Omega} \phi(\rho_\tau(w))\rho_\tau'(w)\psi(x,t,w)\diff{S}\notag\\
\leq&\varepsilon\int_{ \partial\Omega} \widetilde{\Phi}(\phi(\rho_\tau(w)))\diff{S}+C(\varepsilon)\int_{ \partial\Omega} \Phi(|d|)\diff{S}
-\int_{ \partial\Omega} \phi(\rho_\tau(w))\rho_\tau'(w)\psi(x,t,w)\diff{S}\notag\\
\leq&\varepsilon\delta_1\int_{ \partial\Omega} \Phi(\rho_\tau(w))\diff{S}+C(\varepsilon)\int_{ \partial\Omega} \Phi(|d|)\diff{S}
-\int_{ \partial\Omega} \phi(\rho_\tau(w))\rho_\tau'(w)\psi(x,t,w)\diff{S}\notag\\
\leq &\frac{\varepsilon\delta_1}{1+\delta_0}\int_{ \partial\Omega} \phi(\rho_\tau(w))\rho_\tau(w)\diff{S}+C(\varepsilon)\int_{ \partial\Omega} \Phi(|d|)\diff{S}
\notag\\
&-\int_{ \partial\Omega} \phi(\rho_\tau(w))\rho_\tau'(w)\psi(x,t,w)\diff{S}\notag\\
\leq & \frac{\varepsilon\delta_1}{1+\delta_0}\int_{ \partial\Omega} \phi(\rho_\tau(w))\rho_\tau'(w)w\diff{S}+\frac{3\varepsilon\delta_1}{8(1+\delta_0)}\tau\int_{ \partial\Omega} \phi(\rho_\tau(w))\diff{S}\notag\\
&+C(\varepsilon)\int_{ \partial\Omega} \Phi(|d|)\diff{S}
 -\int_{ \partial\Omega} \phi(\rho_\tau(w))\rho_\tau'(w)\psi(x,t,w)\diff{S}\notag\\
\leq & \varepsilon\delta_1\int_{ \partial\Omega} \Phi(\rho_\tau(w))\diff{S}+C(\varepsilon)\int_{ \partial\Omega} \Phi(|d|)\diff{S}\notag\\
&+\frac{3\varepsilon\delta_1}{8(1+\delta_0)}\tau\int_{ \partial\Omega} \phi(\rho_\tau(w))\diff{S} \notag\\
&-\int_{ \partial\Omega} \phi(\rho_\tau(w))\rho_\tau'(w)\psi(x,t,w)\diff{S}.
\end{align}
By \textbf{(A3-2')}, Lemma \ref{property of phi} (iii) and \eqref{property of rho}, it follows that
\begin{align}\label{1602}
-\int_{ \partial\Omega} \phi(\rho_\tau(w))\rho_\tau'(w)\psi(x,t,w)\diff{S}
=&\int_{ \partial\Omega\cap \{w\neq 0\}}\phi(\rho_\tau(w))\rho_\tau'(w)w\cdot \frac{-\psi(x,t,w)}{w} \diff{S}\notag\\
\leq &\int_{ \partial\Omega\cap \{w\neq 0\}} \!\!\!\phi(\rho_\tau(w))\rho_\tau'(w)w \bigg(\!\!\min\limits_{i\in\{0,1\}}\bigg\{\frac{\bm{b}\!\cdot\! \nu }{1+\delta_i}\bigg\}-\underline{\psi_0}\!\bigg)\!\diff{S}\notag\\
= &\int_{\partial\Omega\cap \{\bm{b}\cdot \nu\geq 0 \}}\phi(\rho_\tau(w))\rho_\tau'(w)w \min\limits_{i\in\{0,1\}}\bigg\{\frac{\bm{b}\cdot \nu }{1+\delta_i}\bigg\}\diff{S}\notag\\
&+\int_{\partial\Omega\cap \{\bm{b}\cdot \nu< 0 \}}\phi(\rho_\tau(w))\rho_\tau'(w)w \min\limits_{i\in\{0,1\}}\bigg\{\frac{\bm{b}\cdot \nu }{1+\delta_i}\bigg\}\diff{S}\notag\\
&- \underline{\psi_0}\int_{ \partial\Omega} \phi(\rho_\tau(w))\rho_\tau'(w)w \diff{S}\notag\\
\leq &\int_{\partial\Omega\cap \{\bm{b}\cdot \nu\geq 0\}}\frac{\phi(\rho_\tau(w))\rho_\tau(w)}{1+\delta_1}{\bm{b}\cdot \nu }\diff{S}+\int_{\partial\Omega\cap \{\bm{b}\cdot \nu< 0 \}}\frac{\phi(\rho_\tau(w))\rho_\tau(w)}{1+\delta_0} {\bm{b}\cdot \nu }\diff{S}\notag\\
&-\frac{3}{8} \tau\int_{\partial\Omega\cap \{\bm{b}\cdot \nu< 0 \}}\phi(\rho_\tau(w)) \frac{\bm{b}\cdot \nu }{1+\delta_0}\diff{S}- \underline{\psi_0}\int_{ \partial\Omega} \phi(\rho_\tau(w))\rho_\tau(w) \diff{S}
\notag\\
&+ \frac{3}{8}\underline{\psi_0}\tau\int_{ \partial\Omega} \phi(\rho_\tau(w))\diff{S}
\notag\\
\leq &\int_{\partial\Omega}\Phi(\rho_\tau(w)) {\bm{b}\cdot \nu }\diff{S}- (1+\delta_0)\underline{\psi_0}\int_{ \partial\Omega} \Phi(\rho_\tau(w)) \diff{S}\notag\\
&+ \frac{3}{8} \tau\int_{\partial\Omega} \bigg(\frac{|\bm{b}|}{1+\delta_0}+\underline{\psi_0}\bigg)\phi(\rho_\tau(w))\diff{S}.
\end{align}
Similarly, we have
\begin{align}\label{1507'}
&-\int_{\Omega}c\phi(\rho_\tau(w))w\rho_\tau'(w)\diff{x}
-\int_{\Omega}h(x,t,w,\nabla w)\phi(\rho_\tau(w))\rho_\tau'(w) \diff{x}\notag\\
\leq &-\int_{\Omega}(c-\mu)\phi(\rho_\tau(w))w\rho_\tau'(w)\diff{x}
+\frac{3}{8}\tau \int_{\Omega}(|c|+|\mu|)\phi(\rho_\tau(w))\diff{x}\notag\\
\leq &-\underline{c}\int_{\Omega} \phi(\rho_\tau(w))w\rho_\tau'(w)\diff{x}-\int_{\Omega}
\max\limits_{i\in\{0,1\}}\bigg\{\frac{\divi \ \textbf{b} }{1+\delta_i}\bigg\}\phi(\rho_\tau(w))w\rho_\tau'(w)\diff{x}+\frac{3}{8}\tau \int_{\Omega}(|c|+|\mu|)\phi(\rho_\tau(w))\diff{x}\notag\\
\leq &-\underline{c}(1+\delta_0)\int_{\Omega} \Phi(\rho_\tau(w))\diff{x}-\int_{\Omega}
\Phi(\rho_\tau(w)) \divi \ \textbf{b} \diff{x}+\frac{3}{8}\tau \int_{\Omega}\bigg(\frac{|\divi \ \textbf{b}|}{1+\delta_0}+|c|+|\mu|\bigg)\phi(\rho_\tau(w))\diff{x}.
\end{align}
Note that $\phi'(\rho_\tau(w))(\rho_\tau'(w))^2+\phi(\rho_\tau(w))\rho_\tau''(w)\geq 0$ and
 \begin{align}\label{1511}
\int_{\Omega}f\phi(\rho_\tau(w))\rho_\tau'(w) \diff{x}
\leq & \int_{\Omega}|f|\phi(\rho_\tau(w)) \diff{x}\notag\\
\leq &C(\varepsilon) \int_{\Omega}\Phi(|f|)  \diff{x}+\varepsilon\int_{\Omega}\widetilde{\Phi}(\phi(\rho_\tau(w) ))  \diff{x}\notag\\
\leq &C(\varepsilon) \int_{\Omega}\Phi(|f|)  \diff{x}+\varepsilon\int_{\Omega}\Phi(\rho_\tau(w) )  \diff{x}.
\end{align}
 From \eqref{1500'} to \eqref{1511}, we obtain
\begin{align*}
\frac{\diff{}}{\diff{t}}\int_{\Omega}\Phi(\rho_\tau(w))\diff{x}
\leq-(\underline{c}(1+\delta_0)-\varepsilon)\int_{\Omega} \Phi(\rho_\tau(w))\diff{x}+A(t)+B(t)\tau,
\end{align*}
where we chose $ \varepsilon>0$ small enough such that $ \varepsilon<\min\{ \frac{1+\delta_0}{\delta_1}\underline{\psi_0},\underline{c}(1+\delta_0)\}$, and
\begin{align*}
A(t)=&C(\varepsilon)\int_{ \partial\Omega} \Phi(|d|)\diff{S}+C(\varepsilon) \int_{\Omega}\Phi(|f|)  \diff{x},\\
B(t)=&\frac{3}{8} \int_{\partial\Omega} \bigg(\frac{|\bm{b}|+\delta_1}{1+\delta_0}+\underline{\psi_0}\bigg)\phi(\rho_\tau(w))\diff{S}
+\frac{3}{8} \int_{\Omega}\bigg(\frac{|\divi \ \textbf{b}|}{1+\delta_0}+|c|+|\mu|\bigg)\phi(\rho_\tau(w))\diff{x}.
\end{align*}
The rest part of the proof can be proceeded in the same way as in \eqref{041801}, and \eqref{041802'} can be obtained.

Finally, we show that System \eqref{LPE1'} with the Robin or Neumann boundary condition \eqref{Robin} is exponential $L^q$-ISS in $L^\Phi$-norm w.r.t. boundary and in-domain disturbances $d$ and $f$ from  $L^q_{loc}(\mathbb{R}_+;K_{\Phi}(\partial \Omega))$- and $L^q_{loc}(\mathbb{R}_+;K_{\Phi}(\Omega))$-classes respectively.
Indeed, by Lemma \ref{property of inverse of phi'} (i), we have
\begin{align*}
         \| w(\cdot,T)\|_{L^{\Phi}(\Omega)}\leq & \bigg(\frac{1+\delta_1}{1+\delta_0}\int_{\Omega}\Phi(|w(\cdot,T)|)\diff{x}\bigg)^{\frac{1}{1+\delta_0}}+
\bigg(\frac{1+\delta_1}{1+\delta_0}\int_{\Omega}\Phi(|w(\cdot,T)|)\diff{x}\bigg)^{\frac{1}{1+\delta_1}},\\
    \int_{\Omega}\Phi(|w^0|)\diff{x}   \leq & \frac{1+\delta_1}{1+\delta_0}\big(\| w^0\|_{L^{\Phi}(\Omega)}^{1+\delta_0}+\| w^0\|_{L^{\Phi}(\Omega)}^{1+\delta_1}\big).
\end{align*}
Then we get by \eqref{041802'} and the inequality $(\xi_1+\xi_2)^\epsilon\leq \xi_1^\epsilon+\xi_2^\epsilon\  (\forall \xi_1,\xi_2\in \mathbb{R}_{+},\epsilon\in [0,1])$ that
\begin{align*}
  \| w(\cdot,T)\|_{L^{\Phi}(\Omega)}
  &\leq\!\e^{-\frac{\lambda}{1+\delta_1}T}\!\sum_{i=0}^1\!\bigg( \! \frac{1\!+\!\delta_1}{1\!+\!\delta_0}\!\bigg)^{\!\frac{2}{1+\delta_i}} \!\!\big(\| w^0\|_{L^{\Phi}(\Omega)}^{1+\delta_0}+\| w^0\|_{L^{\Phi}(\Omega)}^{1+\delta_1}\big)^{ \frac{1}{1+\delta_i}}\notag\\
  &+\!\sum_{i=0}^1\!\!\bigg(\! \frac{1\!+\!\delta_1}{1\!+\!\delta_0}C(\varepsilon)\bigg(\! \frac{1}{\lambda q'}\!\bigg)^{\!\!\frac{1}{q'}} \!\bigg)^{\!\!\frac{1}{1+\delta_i}} \!\! \bigg(\!\!\int_0^T\!\!\!\!\bigg(\!\int_{\partial\Omega}\!\!\!\!\Phi(|d|)  \!\diff{S}\!\bigg)^{\!\!q}\!\!\diff{t}\!\!\bigg)^{\!\!\frac{1}{q(1+\delta_i)}}\notag\\
  &+\!\sum_{i=0}^1\!\!\bigg(\! \frac{1\!+\!\delta_1}{1\!+\!\delta_0}C(\varepsilon)\bigg(\! \frac{1}{\lambda q'}\!\bigg)^{\!\!\frac{1}{q'}} \!\bigg)^{\!\!\frac{1}{1+\delta_i}}\! \! \bigg(\!\!\int_0^T\!\!\!\!\bigg(\!\!\int_{\Omega}\!\!\Phi(|f|)  \!\diff{S}\!\bigg)^{\!\!q}\!\!\diff{t}\!\bigg)^{\!\!\frac{1}{q(1\!+\!\delta_i)}}\!\!,
%
%
\end{align*}
which guarantees the exponential $ L^q$-ISS of System \eqref{LPE1'} in $L^\Phi$-norm w.r.t. boundary and in-domain disturbances from  from  $L^q_{loc}(\mathbb{R}_+;K_{\Phi}(\partial \Omega))$- and $L^q_{loc}(\mathbb{R}_+;K_{\Phi}(\Omega))$-classes respectively.
$\hfill\blacksquare$
\end{pf*}
\begin{pf*}{Proof of Theorem~\ref{main result 2'}:} We consider \eqref{L1} with $\mathscr{B}[w]$ given by \eqref{Dirichlet}. As the proof can be proceeded in the same way as that for Theorem~\ref{main result}, we only prove \eqref{ISS-Dirichlet''} and present the main steps.

We choose $ \int_{\Omega}p\Phi(\rho_\tau(w))\diff{x}$ as an approximation of the non-coercive Lyapunov function candidate $\int_{\Omega}p\Phi(|w|)\diff{x}$. Indeed, taking $p\phi( \rho_\tau(w))\rho_\tau'(w)$ as a test function and by direct computations, we have
\begin{align}\label{16a}
\frac{\diff{}}{\diff{t}}\int_{\Omega}p\Phi(\rho_\tau(w))\diff{x}
=&-\int_{ \partial\Omega} \!a\Phi(\rho_\tau(w))\nabla  p \cdot \bm{\nu}\diff{S}+
\int_{\Omega}\!\Phi(\rho_\tau(w))\divi \ (a\nabla p)\diff{x}
\notag\\
&-
\int_{\Omega} \! a{p}\bigg(\!\phi(\rho_\tau(w))(\rho_\tau'(w))^2+\phi(\rho_\tau(w))\rho_\tau''(w)\!\bigg)|\nabla w|^2\diff{x}\notag\\
&+ \int_{\Omega} \Phi(\rho_\tau(w))\divi \ ( p\bm{b})\diff{x}
-\int_{\Omega}cp\phi( \rho_\tau(w))\rho_\tau'(w)w\diff{x}\notag\\
&-\int_{\Omega}ph(x,t,w,\nabla w)\phi( \rho_\tau(w))\rho_\tau'(w)\diff{x}+\int_{\Omega}pf\phi( \rho_\tau(w))\rho_\tau'(w) \diff{x},
\end{align}
where 
\begin{align}
&\int_{\Omega} \Phi(\rho_\tau(w))\divi \ ( p\bm{b})\diff{x}
=\int_{\Omega} \Phi(\rho_\tau(w))\nabla p\cdot\bm{b}\diff{x}+\int_{\Omega} p\Phi(\rho_\tau(w))  \divi \ \bm{b}\diff{x},\label{0412b'} \\
&-\int_{\Omega}cp\phi(\rho_\tau(w))w\rho_\tau'(w)\diff{x}
-\int_{\Omega}ph(x,t,w,\nabla w)\phi(\rho_\tau(w))\rho_\tau'(w) \diff{x}\notag\\
\leq &-\underline{c}(1+\delta_0)\int_{\Omega} p\Phi(\rho_\tau(w))\diff{x}-\int_{\Omega}
p\Phi(\rho_\tau(w)) \divi \ \textbf{b} \diff{x}+\frac{3}{8}\tau \int_{\Omega}p\bigg(\frac{|\divi \ \textbf{b}|}{1+\delta_0}+|c|+|\mu|\bigg)\phi(\rho_\tau(w))\diff{x}, \label{1507'1} \\
&\int_{ \Omega} p\phi( \rho_\tau(w))\rho_\tau'(w)f\diff{x}\notag\\
\leq& C(\varepsilon) \int_{\Omega}\Phi(|pf|)  \diff{x}+\varepsilon\int_{\Omega}\Phi(\rho_\tau(w) )  \diff{x}\notag\\
\leq &C(\varepsilon) \int_{\Omega}\frac{1+\delta_1}{1+\delta_0}\max\{ |p|^{1+\delta_0},|p|^{1+\delta_1}\}\Phi(|f|)  \diff{x}+\varepsilon\int_{\Omega}\Phi(\rho_\tau(w) )  \diff{x}.\label{0901''}
\end{align}
Note that by the boundary condition and {}{\textbf{(A3'-2)}}, we have $w=\psi^{-1}_2\big(\frac{d}{\psi_1}\big)$ on $\partial \Omega\times (0,T)$, which implies
\begin{align}\label{0412a'}
-\int_{ \partial\Omega} a\Phi(\rho_\tau(w))\nabla  p \cdot \bm{\nu}\diff{S}
\leq \int_{ \partial\Omega} a\cdot\bigg((\Phi\circ\rho_\tau\circ\psi^{-1}_2)\left(\frac{d}{\psi_1}\right)\bigg)\cdot|\nabla  p |\diff{S}.
\end{align}
From \eqref{16a} to \eqref{0412a'}, and noting \textbf{(A3'-1)}, we deduce that
\begin{align}\label{0301'''}
\frac{\diff{}}{\diff{t}}\int_{\Omega}p\Phi(\rho_\tau(w))\diff{x}\leq& -\!\big(\underline{c}(1+\delta_0)-\varepsilon\big)\!\!\int_{\Omega} \!p\Phi(\rho_\tau(w))\diff{x}\!+\!A(t)\!+\!B(t)\tau,
\end{align}
where
\begin{align*}
A(t)=&\int_{ \partial\Omega} a\cdot\bigg((\Phi\circ\rho_\tau\circ\psi^{-1}_2)\left(\frac{d}{\psi_1}\right)\bigg)\cdot|\nabla  p |\diff{S}+C(\varepsilon) \int_{\Omega}\frac{1+\delta_1}{1+\delta_0}\max\{ |p|^{1+\delta_0},|p|^{1+\delta_1}\}\Phi(|f|)  \diff{x},\\
 B(t)=&\frac{3}{8}\int_{\Omega}p\bigg(\frac{|\divi \ \textbf{b}|}{1+\delta_0}+|c|+|\mu|\bigg)\phi(\rho_\tau(w))\diff{x}.
\end{align*}
Applying Gronwall's inequality to \eqref{0301'''} and letting $\tau\rightarrow 0$, we may get \eqref{ISS-Dirichlet''}.
$\hfill\blacksquare$
\end{pf*}
\begin{remark}
With more specific arguments, the techniques presented here are still suitable for $g$ satisfying $|g(x,t,s)| \leq g_0|s| $ for all $(x,t,s)\in {\Omega}\times \mathbb{R}_{+}\times \mathbb{R}$, where $g_0$ is a constant. A typical example of such $g$ is: $g(x,t,s)=C_1s+C_2\ln(1+s^2)$, where $C_1,C_2$ are constants. Moreover, one can also conduct ISS analysis by using the technique of transforming for the case of $\underline{c}\geq 0$ in \textbf{(A2-3')}.
\end{remark}

\section{Further comments}\label{Comments}

\subsection{Relaxation of Condition~\textbf{(A2-3)}}
It should be mentioned that although we restrict ourselves to $c-\divi \ \bm{b}-\mu\geq 0$ (or $c-\mu-\max\left\{\frac{\divi \ \textbf{b} }{1+\delta_0},\frac{\divi \ \textbf{b}}{1+\delta_1}\right\}>0$) in this paper, we can still establish ISS (or iISS) estimates for {certain systems for which} $c-\divi \ \bm{b}-\mu$ (or $c-\mu-\max\{\frac{\divi \ \textbf{b} }{1+\delta_0},\frac{\divi \ \textbf{b}}{1+\delta_1}\}$) may change its sign. For example, for System~\eqref{LPE1'} with the Robin or Neumann boundary condition \eqref{Robin}, we have the following result.
\begin{proposition}\label{main result 3}
Suppose that Assumption 1, 3, 4 hold and conditions \textbf{(A2-1)}, \textbf{(A2-2)} in Assumption 2 are satisfied. If $a,\textbf{b},c,\mu$ are independent of the second variable $t$ and $\bm{m}\equiv 0$, then
system \eqref{LPE1'}  with the Robin or Neumann boundary condition \eqref{Robin} is exponential $L^1$-ISS in $L^1$-norm w.r.t. boundary and in-domain disturbances having the following estimate for all $T>0$ and any constant $\underline{c}> \max\{\max\limits_{\overline{\Omega}}( c-\divi \ \bm{b}-\mu),0\}$:
\begin{align*}
\|w(\cdot,T)\|_{L^1( \Omega)}
\leq  &\e^{\|u\|_{L^{\infty}(\Omega)}}\e^{-\underline{c}T}\|{w^0}\|_{L^1( \Omega)}
+ \e^{\|u\|_{L^{\infty}(\Omega)}}\| {d}\|_{L^1(\partial \Omega \times (0,T))}
+\e^{\|u\|_{L^{\infty}(\Omega)}}\| {f}\|_{L^1( Q_T)},
\end{align*}
where $u\in C^2(\overline{\Omega})$ satisfies the elliptic equation \eqref{elliptic} in Appendix \ref{Proof of Proposition}.
\end{proposition}
 The proof of Proposition \ref{main result 3} is given in Appendix \ref{Proof of Proposition}.

\begin{remark}
In general, it is not an easy task to obtain the boundedness of $u$ for the nonlinear elliptic equation \eqref{elliptic}. However, it is possible to obtain the boundedness of the solution for some special equations. For example,  
let $\Omega=B_1(0)$ and assume that
\begin{align*}
K=1,a(x)\equiv A(|x|), \ \bm{b}=b_0(|x|)x,\divi \ \bm{b}-c+{\mu}+\underline{c}=B(|x|),
\end{align*}
where $A,b_0\in C^{2}(\mathbb{R}_{\geq 0};\mathbb{R}),B\in C^{2}(\mathbb{R}_{\geq 0};\mathbb{R}_{\geq 0})$ are given functions with $A(1)=1$. It is expected to find a radially symmetric solution $u(x)=u(|x|)$ of System~\eqref{elliptic}. Indeed, letting $r=|x|$, it follows that $\nabla u=\frac{u_r}{r}x,\Delta u=u_{rr}+\frac{n-1}{r}u_r$. Let $Q(r)=\frac{1}{A}\left(B-A_r+\frac{(n-1)A}{r}-b_0r\right)$. System~\eqref{elliptic} becomes:
\begin{subequations}\label{elliptic'}
\begin{align}
&u_{rr}=u_r^2+Q(r),\ \ r\in (0,1),\\
&u_r(0)=0, u_r(1)+u(1)=0,
\end{align}
\end{subequations}
which may be solvable for certain $Q(r)$ as shown in Section~\ref{special equation}.
\end{remark}
\subsection{A remark on ISS in weighted $L^1$-norm for parabolic PDEs with Robin boundary conditions}

In \cite{karafyllis2017siam}, the authors considered the following $1$-$D$ linear parabolic equation with $L^\infty$-Dirichlet disturbances:
\begin{subequations}\label{example linear PDE}
\begin{align}
&\!\!\!\!w_t-aw_{xx}+bw_x+cw=f(x,t), (x,t)\in(0,1)\!\times\!\mathbb{R}_+,\label{example linear PDE a}\\
&\!\!\!\!w(0,t)=d_0(t),w(1,t)=d_1(t), t\in \mathbb{R}_{+},\\
&\!\!\!\!w(x,0)=w^0(x), x\in (0,1),
\end{align}
\end{subequations}
where $d_0,d_1,f,w^0$ are given functions, and $a>0,b,c\in \mathbb{R}$ are constants satisfying
\begin{align}
 c+a\pi^2+\frac{b^2}{4a}>0.\label{c}
\end{align}
Under certain regularity and compatibility conditions on $d_0,d_1,f$, and $w^0$, and applying a finite-difference scheme for the approximation of the solution of \eqref{example linear PDE}, {the following ISS estimate in weighted $L^1$-norm is obtained (see \cite[Theorem 2.4]{karafyllis2017siam})}:
\begin{align}\label{weight-Dirichlet}
\int_{0}^1\e^{-\frac{bx}{2a}}\sin(\pi x)|w(x,T)|\diff{x}
\leq&\frac{4a^2\pi}{4a^2\pi^2\!+\!b^2\!+\!4ac}\!\left(\!\!\left(\!\max\limits_{t\in[0,T]}|d_0(t)|\!\right)
 \!+\!\e^{-\frac{b}{2a}}\!\max\limits_{t\in[0,T]}|d_1(t)|\!\right)\notag\\
&\!\!\!+\!\frac{4a}{4a^2\pi^2\!+\!b^2\!+\!4ac}\!\max\limits_{t\in[0,T]}\!\!\left( \! \max\limits_{x\in[0,1]}|f(x,t)|\e^{-\frac{bx}{2a}}\sin (\pi x)\!\!\right)\notag\\
&\!\!\!+ \!\e^{-\!\left(\!a\pi^2+c+\frac{b^2}{4a}\!\right)T}\!\!\!\!
\int_{0}^1\!\!\!\e^{-\frac{bx}{2a}}\sin(\pi x)|w^0(x)|\diff{x},\!\forall T\!>\!0.
\end{align}
Meanwhile, it is claimed that \eqref{weight-Dirichlet} is only applicable to \eqref{example linear PDE a}  with Dirichlet boundary conditions ( see \cite[page 1721]{karafyllis2017siam}).

We would like to point out that if we replace \eqref{c} with the the following condition:
\begin{align*}
 \exists\theta\in \left[0,\frac{\pi}{2}\right) \ \text{s.t.}\ c+a\theta^2+\frac{b^2}{4a}>0,
\end{align*}
a similar ISS estimate in weighted $L^1$-norm is applicable to \eqref{example linear PDE a}  with Robin boundary conditions. For example, let $w$ be the solution of \eqref{example linear PDE a} with the following Robin boundary conditions for $t\geq 0$:
\begin{align*}
w_x(0,t)-Kw(0,t)=d_0(t), w_x(1,t)+Kw(1,t)=&d_1(t),
\end{align*}
where $K>0$ is a sufficiently large constant. 
Let $u=\beta w$ with $\beta(x)=\e^{-\frac{bx}{2a}}\cos(\theta x)>0$. It is easy to check that $u $ satisfies:
\begin{align*}
&u_{t}-au_{xx}+Bu_x+Cu=F(x,t), (x,t)\in(0,1)\times\mathbb{R}_{+},\label{example linear PDE u a}\\
&u_x(0,t)-K_0u(0,t)=D_0(t), t\in\mathbb{R}_{+},\\
&u_x(1,t)+K_1u(1,t)=D_1(t), t\in\mathbb{R}_{+},\\
&u(x,0)=u^0(x),x\in (0,1),
\end{align*}
where 
$B=2a\frac{\beta_x}{\beta}+b,C= c-b\frac{\beta_x}{\beta}+a\left(\frac{\beta_{xx}}{\beta}-2\frac{\beta_x^2}{\beta^2}\right),F=\beta f,K_0=K- \frac{b}{2a} ,K_1=K+\frac{b}{2a}+\theta \tan \theta$, $D_i(t)=\beta(i)d_i(t),i=0,1$, and $u^0=\beta w^0$.

By direct computations, we have $C-\divi \ B=c-b\frac{\beta_x}{\beta}-a\frac{\beta_{xx}}{\beta}=c+a\theta^2+\frac{b^2}{4a}:=\underline{c}>0$.  
 Assume that $K$ is large enough such that $K_i>0,i=0,1$. According to Theorem~\ref{main result}~(i), for any $ T>0$, we have (see \eqref{1803}):
\begin{align*}
\|u(\cdot,T)\|_{L^{\!1}( 0,1)}
\leq  &\e^{-\underline{c}T}\|u^0\|_{L^1( 0,1)}+\frac{1}{\underline{c}}\!\bigg( \! \max\limits_{t\in[0,T]}\!|\beta(0)d_0(t)|\!+\!\max\limits_{t\in[0,T]}\!|\beta(1)d_1(t)|+\max\limits_{(x,t) \in[0,1]\times[0,T]}|\beta (x)f(x,t)|\bigg).
\end{align*}
Therefore, the following ISS estimate in weighted $L^1$-norm holds:
\begin{align*}
\int_{0}^1\e^{-\frac{bx}{2a}}\cos(\theta x) |w(x,T)|\diff{x}
\leq  &\e^{-\left(a\theta^2+c+\frac{b^2}{4a}\right)T}\int_{0}^1\e^{-\frac{bx}{2a}}\cos(\theta x) |w^0(x)|\diff{x}\notag\\
&\!+\!\frac{4a}{4a^2\theta^2\!+\!b^2\!+\!4ac}\left( \! \max\limits_{t\in[0,T]}\!|d_0(t)|\!+\!\e^{-\frac{b}{2a}}\!\max\limits_{t\in[0,T]}|d_1(t)|\!\right)\notag\\
&\!+\!\frac{4a}{4a^2\theta^2\!+\!b^2\!+\!4ac}\!\left(\! \max\limits_{(x,t)\in [0,1]\!\times\![0,T]}\!\e^{-\frac{bx}{2a}}\!\cos(\theta x) |f(x,t)|\!\right).
\end{align*}

Thus the results obtained in this paper can be seen as a complement of  \cite[Theorem 2.4]{karafyllis2017siam}.

\subsection{Remarks on the main results}\label{remarks existence}

Due to the restriction of the special {}{approximations of Lyapunov functions} constructed in this paper, the ISS or iISS estimates for System \eqref{LPE1'} with Dirichlet boundary disturbances \eqref{Dirichlet} involve weighted norms or weighted $E_\Phi$-classes, which can be seen as extensions of
\cite[Theorem 2.4 and Corollary 2.5]{karafyllis2017siam} to nonlinear PDEs over higher dimensional domains. It is worthy noting that it is possible to {assess the} $L^q$-ISS without a weighted norm by constructing some non-coercive Lyapunov functions as in \cite{Jacob:2019} and applying the approximating function $\rho_\tau(\cdot)$ in the same time.


{Furthermore, in} this paper, the data $w^0,d$, and $f$ are supposed to be smooth enough so that the $L^1$-ISS estimates in $L^1$-norm can be established for classical solutions. If the data $w^0,d$, and $f$ belong only to $L^1$-space, the functions $|a\nabla w|$ and $N[w](x,t)$ do not belong to $L^1((0,T)\times \Omega)$ {}{and $L^1((0,T)\times \partial\Omega)$} in general. In this case, it will be an arduous task to prove existence of a weak solution even for the elliptic equations (see e.g., \cite{Benilan:1995,Boccardo:1999,Zheng:2019obstacle,Zheng:2019obstacle2}). Besides, the relaxation on the regularity of the data $w^0,d$, and $f$ {may also significantly influence} the regularity of the solutions, which in turn causes some crucial difficulties in {stability investigation. To overcome this issue,} a possible solution is the use of the notion of \emph{renormalized solution} or \emph{entropy solution}. It should be mentioned that the notion of renormalized solution was first introduced by Di~Perna and Lions \cite{DiPerna:1989} in the study of Boltzmann equation (see also \cite{Lions:1996} for applications to fluid mechanics models). The equivalent notion of entropy solution was developed independently in \cite{Benilan:1995} for nonlinear elliptic problems. It is worthy noting that using the function $\rho_\tau(\cdot)$ introduced in this paper, we can also establish ISS estimates for nonlinear PDEs with boundary or in-domain disturbances from $L^1$-space in the framework of renormalized solutions or entropy solutions, which will be addressed in our future work.


\section{Illustrative examples}\label{Sec: examples}

\subsection{Generalized Burgers' equation with nonlinear Robin boundary disturbances}
We consider the following $n$-dimensional generalized Burgers' equation with Robin boundary disturbances over $B_R(0)\times \mathbb{R}_{+}$:
\begin{subequations}\label{Burgers}
\begin{align}
 &\!\!w_t\!-\!\Delta w\!+\!\lambda w\!+\!\sum_{i=1}^nm_iw w_{x_i} \! =\!f, (x,t)\!\in\! B_R(0)\!\times\! \mathbb{R}_+,\\
&\!\!\frac{\partial w}{\partial\bm{\nu}}+K(w+w^3)=d, (x,t)\in  \partial B_R(0) \times \mathbb{R}_+,\\
&\!\!w(x,0)=w^0(x),x\in  B_R(0),
\end{align}
\end{subequations}
where $\lambda\geq0,K>0$, and $m_i (i=1,2,...,n)$ are constants.

For the above system, we set
\begin{align*}
 &a(x,t)\equiv 1, \bm{b}(x,t)\equiv\bm{0},c(x,t)\equiv \lambda,\notag\\
&h(x,t,w,\nabla w)\equiv\mu(x,t)\equiv0,g(x,t,w)\equiv \frac{1}{2}w^2,\\
 &\psi(x,t,w)\equiv K(w+w^3),\bm{\nu}(x)\equiv\frac{x}{R},c-\mu-\divi \ \bm{b}\equiv\lambda.
\end{align*}
Moreover, assume that \textbf{(A4)} is satisfied.

$\bm{L^1}$\textbf{-estimates for}  $\bm{\lambda\geq 0}.$
 If we set $\bm{m}(x,t)\equiv(1,1,...,1)$, we have $g(x,t,w)w\cdot\divi \ \bm{m}\equiv0$. Note that
\begin{align*}
|s|\leq 1+s^2\ \text{and}\ \frac{1}{R}\sum\limits_{i=1}^nx_i \leq \frac{1}{R}nR=n, \forall s\in \mathbb{R},x\in\partial B_R.
 \end{align*}
 If $K\geq \frac{n}{2}$, it follows that
\begin{align*}
\psi_0(x,t):=&\inf\limits_{s\in\mathbb{R}\setminus\{0\}}\frac{\psi(x,t,s)
+s\bm{b}\cdot\bm{\nu}+g(x,t,s)\bm{m}\cdot\bm{\nu}}{s}
=\inf\limits_{s\in\mathbb{R}\setminus\{0\}}\frac{K(s+s^3)
+ \frac{1}{2}s^2\frac{1}{R}\sum\limits_{i=2}^nx_i}{s}
\geq
 K-\frac{n}{2}\notag\\
 \geq& 0,\forall (x,t,s)\in \partial{B}_R\times \mathbb{R}_{+}\times \mathbb{R}.
\end{align*}
According to Theorem~\ref{main result}~(i), we know that System~\eqref{Burgers} has the $L^1$-estimate \eqref{042003}  with $\underline{c}=\lambda$. Particularly, for $\lambda>0$, System~\eqref{Burgers} is exponential $L^q$-ISS in $L^1$-norm for any $q\in [1,+\infty]$. Furthermore, by Theorem~\ref{main result Phi}~(i), System~\eqref{Burgers} has the $L^1$-estimate \eqref{4.3} with $\underline{c}=\lambda$. Moreover, for $\lambda>0$, System~\eqref{Burgers} is exponential $L^\Phi$-ISS in $L^1$-norm for any $\Phi$ defined in Section~\ref{ISS Orlicz}.

$\bm{L^q}$\textbf{-ISS and} $\bm{L^\Phi}$\textbf{-ISS in} $\bm{L^1}$\textbf{-norm for} $\bm{\lambda\geq 0}.$
If we set $\bm{m}(x,t)\equiv(0,1,...,1)$ and $n\geq 2$, we have
\begin{align*}
g(x,t,w)w\!\cdot\!\divi \ \bm{m}\equiv\bm{m}_1\equiv
&\frac{b_j+a_{x_j}}{a}\equiv0,\forall j\!\in\!\{1,2,...,n\}.
\end{align*}
Similarly, if $K>\frac{n-1}{2}$,
 it follows that
\begin{align*}
\psi_0(x,t):=&\inf\limits_{s\in\mathbb{R}\setminus\{0\}}\frac{\psi(x,t,s)
+s\bm{b}\cdot\bm{\nu}+g(x,t,s)\bm{m}\cdot\bm{\nu}}{s}
\notag\\
&=\inf\limits_{s\in\mathbb{R}\setminus\{0\}}\frac{K(s+s^3)
+ \frac{1}{2}s^2\frac{1}{R}\sum\limits_{i=2}^nx_i}{s}\notag\\
&
\geq K-\frac{n-1}{2}
:=\underline{\psi_0}>0.
\end{align*}
According to Theorem~\ref{main result}~(ii), we know that System~\eqref{Burgers} is exponential $L^q$-ISS for any $q\in [1,+\infty]$, with $\widehat{\underline{c}}=\frac{l\e^{-l\mathbbm{d}}}{k +\e^{l\mathbbm{d}}}\big(l-\frac{2}{k}l\e^{l\mathbbm{d}}\big)$ and $C_{k,l}=\frac{k+\e^{l\mathbbm{d} }}{k+\e^{-l\mathbbm{d} }}$, where $l$ and $k$ can be any constants satisfying $l>0$ and $k>\max\{ 2\e^{l\mathbbm{d}},\frac{ l\e^{l\mathbbm{d}}}{\underline{\psi_0}}\}$. Furthermore, by Theorem~\ref{main result Phi}~(ii), System~\eqref{Burgers} is exponential $L^\Phi$-ISS in $L^1$-norm for any $\Phi$ defined in Section~\ref{ISS Orlicz}.

$\bm{L^q}$\textbf{-ISS in}  $\bm{L^\Phi}$\textbf{-norm for} $\bm{\lambda> 0}.$
If we set $\bm{m}(x,t)\equiv \bm{0}$ and $K>0$, it follows that
\begin{align*}
\psi_0(x,t):=&\inf\limits_{s\in\mathbb{R}\setminus\{0\}}\frac{\psi(x,t,s)
+s\min\{\frac{\bm{b}\cdot\bm{\nu}}{1+\delta_0},\frac{\bm{b}\cdot\bm{\nu}}{1+\delta_1}\}}{s}
=\inf\limits_{s\in\mathbb{R}\setminus\{0\}}\frac{K(s+s^3)}{s}
=K>0,\forall (x,t,s)\in \partial{B}_R\times \mathbb{R}_{+}\times \mathbb{R}.
\end{align*}
According to Theorem~\ref{main result'}, we know that System~\eqref{Burgers} is exponential $L^q$-ISS in $L^\Phi$-norm for any $q\in [1,+\infty]$ and $\Phi$ defined in Section~\ref{ISS Orlicz}.

\subsection{Generalized Ginzburg-Landau equation with nonlinear Dirichlet boundary disturbances}
We consider the following generalized Ginzburg-Landau equation
with variable coefficients and nonlinear Dirichlet boundary disturbances:
\begin{subequations}\label{Chafee-Infante}
\begin{align}
 &w_t\!-\!\divi  \ (a\nabla\! w)
\!+\!c_1w\!+\!c_2 w^3\!+\!c_3 w^5  \!=\!f,(x,t)\in \Omega\times \mathbb{R}_+,\\
&w+w^3=d, (x,t)\in  \partial \Omega \times \mathbb{R}_+,\\
&w(x,0)=w^0(x),  (x,t)\in  \Omega,
\end{align}
\end{subequations}
where $c_1 \in C^2(\mathbb{R};\mathbb{R}_{+}),c_2,c_3\in C^2(\mathbb{R};\mathbb{R}_{\geq 0})$ are given functions.

For the above system, we set
\begin{align*}
 &\bm{b}(x,t)\equiv\bm{0},c(x,t)\equiv c_1(x,t),\notag\\
&
  h(x,t,w,\nabla w)\equiv c_2 (x,t)w^3+c_3(x,t) w^5 ,
 \bm{m}(x,t)\equiv\bm{0}, \notag\\
&g(x,t,w)\equiv 0,
 \psi(x,t,w)\equiv\psi_1(x,t)\psi_2(w)\equiv w+w^3.
\end{align*}
Moreover, assume that \textbf{(A4)} is satisfied.

Note that
$
  -h(x,t,w,\nabla w)w  \leq0$ with $\mu\equiv 0$, $\psi_1(x,t)\equiv1,\psi_2(s)\equiv s+s^3
 $ with $\psi^{-1}_2(s)$ 
 satisfying
 $|\psi^{-1}_2(s)|\leq |s|$ for all $s\in \mathbb{R}$, and $\psi_0(x,t):=\inf\limits_{s\in\mathbb{R}\setminus\{0\}}\frac{\psi(x,t,s)
+s\min\{\frac{\bm{b}\cdot\bm{\nu}}{1+\delta_0},\frac{\bm{b}\cdot\bm{\nu}}{1+\delta_1}\}}{s}\equiv 1$.

\textbf{ISS and iISS estimates for} $\bm{c_1>0}.$
Assume that there exists a positive constant $\underline{c}>0$ such that
\begin{align*}
c_1\geq \underline{c},\forall (x,t)\in\overline{\Omega}\times \mathbb{R}_{\geq 0}.
\end{align*}
For any weighting function $p\in C^2(\Omega)\cap C^1(\overline{\Omega})$ and any constant $p_0\geq0$ satisfying
\begin{align*}
 \divi \ (a\nabla p)&\leq -p_0, (x,t)\in  \Omega\times \mathbb{R}_{+},\\
p&=0, x\in  \partial \Omega,
\end{align*}
the following statements hold true:
\begin{enumerate} [(i)]
 \item according to Theorem~\ref{main result 2}~(i), System~\eqref{Chafee-Infante} is $L^q$-ISS in weighted $L^1$-norm for any $q\in [1,+\infty]$;
 \item according Theorem~\ref{main result 2 Phi}~(i), System~\eqref{Chafee-Infante}  is $L^\Phi$-ISS in weighted $L^1$-norm for any $\Phi$ defined in Section \ref{ISS Orlicz};
 \item  according Theorem~\ref{main result 2'}, System~\eqref{Chafee-Infante} is exponential $L^q$-ISS in weighted $K_{\Phi}$-class w.r.t. boundary and in-domain disturbances $d$ and $f$ from $L^q_{loc}(\mathbb{R}_+;K_{\Phi}(\partial \Omega))$- and $L^q_{loc}(\mathbb{R}_+;K_{\Phi}(\Omega))$-classes respectively, where $\Phi$ is defined in Section \ref{ISS Orlicz}.
\end{enumerate}

\subsection{A special nonlinear parabolic equation with Robin boundary disturbances}\label{special equation}
We consider the following $1$-$D$ nonlinear equation:
\begin{subequations}\label{2002}
\begin{align}
 &w_t- w_{xx}+3xw_x+\bigg(3|x|^2+\frac{5}{2}\bigg)w + \bigg(|x|^2-\frac{3}{\ln 2}|x|\ln(1+|x|)\bigg)\frac{w}{1+w^2}
=f, (x,t)\in (-1,1)\times \mathbb{R}_+,\\
&w_x(1,t)+w(1,t)
=d(1,t),
w_x(0,t)-w(0,t)
=-d(0,t), t\in  \mathbb{R}_+,\\
&w(x,0)=w^0(x), x\in  (-1,1),
\end{align}
\end{subequations}
For the above system, we set
\begin{align*}
 &a(x,t)\equiv 1, b(x,t)\equiv3x,c(x,t)\equiv 3|x|^2+\frac{5}{2},
  \notag\\
&h(x,t,w,w_x)\equiv \bigg(|x|^2-\frac{3}{\ln 2}|x|\ln(1+|x|)\bigg)\frac{w}{1+w^2},\\
 &m(x,t)\equiv0,g(x,t,w)\equiv 0,\psi(x,t,w)\equiv w.
\end{align*}
Assume that \textbf{(A4)} is satisfied. Note that $h\in C^2([-1,1];\mathbb{R})$ satisfying $ -h(x,t,w,w_x)w=-\big(|x|^2-\frac{3}{\ln 2}|x|\ln(1+|x|)\big)\frac{w^2}{1+w^2}\leq (-|x|^2+3|x|)w^2$, we can choose $\mu(x,t):=-|x|^2+3|x|$.

Obviously, $c-b_x-\mu =4|x|^2-3|x|-\frac{1}{2}$ changes its sign on $[-1,1]$. Therefore, Theorem~\ref{main result}~(i) and Theorem~\ref{main result}~(ii) cannot be applied for this case. However, there exists a constant $\underline{c}=\frac{3}{2}$ such that
$\underline{c}> \frac{17}{16}=\max\limits_{[0,1]}| c-b_x-\mu|.$
Thus, we consider \eqref{elliptic} with $n=1,\Omega=(-1,1),K\equiv 1$ and its radially
symmetric solution.

Note that
$B-A_r+\frac{(n-1)A}{r}-b_0r=B-b_0r=2-4r^2,$
where $
 A(|x|)\equiv 1,B(|x|)\equiv b_x-c+{\mu}+\underline{c}.
$
It is easy to see that $u(r)=r^2-3$ is the solution of the following equation:
 \begin{align*}
&u_{rr}=u_r^2+2-4r^2,\ \ r\in (0,1),\\
&u_r(0)=0, u_r(1)+u(1)=0.
\end{align*}
Thus, the lower and upper bounds of $u(x)=|x|^2-3$ is given by
\begin{align}\label{u}
-3\leq u\leq -2.
\end{align}
According to Proposition \ref{main result 3} and \eqref{u}, we deduce that System \eqref{2002} is $L^1$-ISS in $L^1$-norm having the estimate for all $T>0$:
\begin{align*}
\|w(\cdot,T)\|_{L^1( 0,1)}
\leq  &\e^{3}\!\Big(\!\e^{-\frac{3}{2}T}\!\|{w^0}\|_{L^1( 0,1)}
\!+ \!\| {d}(0,\cdot)\|_{L^1(0,T)}\!+\!\| {d}(1,\cdot)\|_{L^1(0,T)}+\| {f}\|_{L^1( (0,1)\times (0,T))}\Big).
\end{align*}

%
%

\section{Conclusion}\label{Sec: Conclusion}
{This paper introduced a convex function that can be used for the construction of approximations of coercive or non-coercive ISS Lyapunov functions.} Based on these approximating functions, ISS and iISS estimates in various norms, in particular $L^1$-norm, for a class of nonlinear parabolic PDEs over higher dimensional domains with different types of boundary disturbances from different spaces or Orlicz classes have been established {via solely the Lyapunov method}. Several examples are provided to illustrate the effectiveness of the developed approach for ISS analysis of nonlinear PDEs with in-domain and different boundary disturbances.

\appendix
\section{Well-posedness for a class of quasi-linear parabolic equations} \label{Appendix I}

Let $\Omega$ be an open bounded domain in $\mathbb{R}^n (n\geq 1)$  with a $C^2$-boundary $\partial \Omega$. We consider the following parabolic equation:
\begin{subequations}\label{quasilinear PDE Dirichlet}
\begin{align}
 &u_t\!-\!\divi  \bm{A}(x,t,u,\!\nabla u)\!+\!B(x,t,u,\!\nabla  u)=0,  (x,t)\in\Omega\!\times\! \mathbb{R}_+,\label{quasilinear PDE Dirichlet a}\\
&\mathscr{B}[u](x,t)=0, (x,t)\in \partial \Omega \times \mathbb{R}_{+},\\
&u(x,0)=u^0,  x\in\Omega,
\end{align}
\end{subequations}
where $\bm{A}(x,t,u,\nabla u)$$=$$(A_1(x,t,u,\nabla u), A_2(x,t,u,\nabla u)$,..., $A_n(x,t,u,\nabla u))$, $A_i\in C^{2}(\overline{\Omega}\times \mathbb{R}_{\geq 0}\times \mathbb{R}\times \mathbb{R}^n; \mathbb{R}), i=1,2,...,n$, $B \in C^{2}(\overline{\Omega}\times \mathbb{R}_{\geq 0}\times \mathbb{R} \times \mathbb{R}^n; \mathbb{R})$  and $u^0\in C^{2}( \overline{\Omega}; \mathbb{R})$ are given functions, $\mathscr{B}[u](x,t)=0$ represents the boundary condition, which is given by 
 \begin{align}\label{Dirichlet condition}
 \!\!\mathscr{B}[u](x,t):= u(x,t)-\Psi(x,t)=0, (x,t)\in \partial \Omega \times \mathbb{R}_{+},
\end{align}
 where $\Psi\in C^{2}( \overline{\Omega}\times \mathbb{R}_{\geq 0}; \mathbb{R})$,
  or 
  \begin{align}\label{Oblique condition}
 \mathscr{B}[u](x,t):=&\sum\limits_{i=1}^nA_i(x,t,u,\nabla u)\cos (\bm{n},x_i)+\Psi(x,t,u))
=0,(x,t)\in \partial \Omega \times \mathbb{R}_{+},
\end{align}
where $\Psi\in C^{2}( \overline{\Omega}\times \mathbb{R}_{\geq 0}\times \mathbb{R}; \mathbb{R})$ and $\bm{n}=\bm{n}(x)$ is {a vector outward normal} to $\partial \Omega$.

Let $\widetilde{B}(x,t,u,\nabla u):=B(x,t,u,\nabla u)-\sum\limits_{i=1}^n\frac{\partial A_i}{\partial u}u_{x_i}-\sum\limits_{i=1}^n\frac{\partial A_i}{\partial x_i}$. We impose the following assumptions:\\
      \textbf{(A5-1)} There exist $\mu_1,\mu_2\in C(\overline{\Omega}\times \mathbb{R}_{\geq 0}\times \mathbb{R}; \mathbb{R}_{+})$ such that
$ \mu_0(x,t,u) \xi^2\leq \sum\limits_{i,j=1}^n\frac{\partial A_i(x,t,u,\eta)}{\partial \eta_j}\xi_i\xi_j\leq \mu_1(x,t,u) \xi^2 $  for all $(x,t,u,\eta)\in\overline{\Omega}\times\mathbb{R}_{\geq 0}\times \mathbb{R}\times \mathbb{R}^n$ and all $\xi \in \mathbb{R}^n$.\\
    \textbf{(A5-2)} There exist $b_1\in C( \overline{\Omega}\times\mathbb{R}_{\geq 0}; \mathbb{R}),b_2\in C( \overline{\Omega}\times\mathbb{R}_{\geq 0}; \mathbb{R}_{\geq 0} )$ such that
$-\widetilde{B}(x,t,u,0)u\leq b_1(x,t)+b_2(x,t)H(|u|)|u|$ for all $ (x,t,u,\eta)\in\overline{\Omega}\times\mathbb{R}_{\geq 0}\times \mathbb{R}\times \mathbb{R}^n$, where $H\in C(\mathbb{R}_{\geq 0};\mathbb{R}_+)$ satisfying $ \int_{0}^{+\infty}\frac{1}{H(s)}\diff{S}=\infty$.\\
  \textbf{(A5-3)} There exists $\mu_2\in C(\overline{\Omega}\times \mathbb{R}_{\geq 0}\times \mathbb{R}; \mathbb{R}_{\geq 0})$ such that
\begin{align*}
\sum\limits_{i=1}^n\left( |A_i|+\left|\frac{\partial A_i}{\partial u}\right|\right)(1+|\eta|)+\sum\limits_{i,j=1}^n\left|\frac{\partial A_i}{\partial x_j}\right|+|b|
\leq &\mu_2(x,t,u)(1+|\eta|)^2,\forall (x,t,u,\eta)\in\overline{\Omega}\times\mathbb{R}_{\geq 0}\times \mathbb{R}\times \mathbb{R}^n.
 \end{align*}
\textbf{(A5-4)} $\mathscr{B} [u^0](x,t)=0$ for all $(x,t)\in \Omega \times \{0\}.$\\
\textbf{(A5-5)}
There exist $\psi_1,\psi_2\in C(\overline{\partial \Omega}\times \mathbb{R}_{\geq 0})$ such that
$-u\Psi(x,t,u)\leq \psi_1(x,t)+\psi_2(x,t)u^2$ for all $(x,t)\in \partial\Omega \times \{0\}.$




%

\begin{theorem}\label{well-posedness result} The following statements hold true:
\begin{enumerate} [(i)]
\item
Under Assumption  \textbf{(A5-1)},  \textbf{(A5-2)}, \textbf{(A5-3)} and \textbf{(A5-4)}, 
\eqref{quasilinear PDE Dirichlet} with the boundary condition \eqref{Dirichlet condition} admits a unique solution $ u\in C^{2,1}(Q_T)\cap C(\overline{Q}_T)$ for any
$T > 0$.
\item
Under Assumption  \textbf{(A5-1)},  \textbf{(A5-2)}, \textbf{(A5-3)}, \textbf{(A5-4)}, 
and \textbf{(A5-5)}, \eqref{quasilinear PDE Dirichlet} with the boundary condition \eqref{Oblique condition} admits a unique solution $ u\in C^{2,1}(Q_T)\cap C(\overline{Q}_T)$ for any
$T > 0$.
\end{enumerate}
\end{theorem}
\begin{pf*}{Proof:} Note that \eqref{quasilinear PDE Dirichlet a} can be written as
$$u_t - \sum\limits_{i,j=1}^n\frac{\partial A_i(x,t,u,\eta)}{\partial \eta_j}u_{x_ix_j} + \widetilde{B}(x,t,u,\nabla  u)=0.$$
Due to the smoothness assumptions on $\mu_0,\mu_1,\mu_2$ and $b_1,b_2,\psi_1,\psi_2$, for any $M>0$ and $T>0$,  $\mu_0,\mu_1,\mu_2$, and $b_1,b_2,\psi_1,\psi_2$ are bounded on $\overline{\Omega}\times [0,T]\times [-M,M]$ and $\overline{\Omega}\times [0,T]$. Particularly, there exist positive constants $\underline{\mu_0},\overline{\mu_1},\overline{\mu_2}$ and $\overline{b_1},\overline{b_2},\overline{\psi_1},\overline{\psi_2}$
such that $\mu_0\geq \underline{\mu_0}, \mu_1\leq \overline{\mu_1},\mu_2\leq \overline{\mu_2}$ on $\overline{\Omega}\times [0,T]\times [-M,M]$ and $b_1\leq \overline{b_1},0\leq b_2\leq \overline{b_2}, \psi_1\leq \overline{\psi_1},\psi_2\leq \overline{\psi_2}$   on $\overline{\Omega}\times [0,T]$, respectively. Therefore, every condition of \cite[Theorem 6.2 Chapter V]{Ladyzhenskaya:1968} is satisfied. Finally, the existence of a unique solution $ u\in C^{2,1}(Q_T)\cap C(\overline{Q}_T)$ to \eqref{quasilinear PDE Dirichlet} with the boundary condition \eqref{Dirichlet condition} is guaranteed by \cite[Theorem 6.1, Chapter V]{Ladyzhenskaya:1968}. The existence and uniqueness of a solution $u\in C^{2,1}(Q_T)\cap C(\overline{Q}_T)$ to \eqref{quasilinear PDE Dirichlet} with the boundary condition \eqref{Oblique condition} can be proceeded exactly as in the proof of \cite[Theorem 7.4, Chapter V]{Ladyzhenskaya:1968} (see also the comments in the last paragraph of \cite[\S 7, Chapter V]{Ladyzhenskaya:1968}).
$\hfill\blacksquare$
\end{pf*}

\section{Proof of Proposition~\ref{main result 3}} \label{Proof of Proposition}
For any constant $\underline{c}>\max\{\max\limits_{\overline{\Omega}}( c-\divi \ \bm{b}-\mu),0\}$, let $u$ be the unique solution of the following elliptic equation:
\begin{subequations}\label{elliptic}
\begin{align}
&-\divi \ (a\nabla u)- \bm{b}\cdot\nabla u
=-a|\nabla u|^2- \divi \ \bm{b}+c-{\mu}-\underline{c},\ \ x\in \Omega,\\
&a \frac{\partial u}{\partial\bm{\nu}}+Ku=0,\ \ x\in \partial \Omega,
\end{align}
\end{subequations}
where $K>0$ is a constant.

Note that by the classical theory of seconder order elliptic PDEs (see, e.g., \cite[Chapter 10]{Ladyzhenskaya:1968b}), $u$ is bounded in $C^{2}(\overline{\Omega})$. Moreover, by the maximum principle 
we always have
$u\leq 0$ for $ x\in \overline{\Omega}.$
Define the following quantities:
\begin{align*}
&v=\e^{-u}w,\widehat{f}=\e^{-u}f,\widehat{d}=\e^{-u}d,\widehat{w}^0=\e^{-u}w^0,
\notag\\
&\widehat{\bm{b}}=\bm{b}-2a\nabla u, \widehat{c}=-  \divi \ (a\nabla u)-a|\nabla u|^2+\bm{b}\cdot\nabla u+c,\\
&\widehat{h}(x,t,s,\eta)=\e^{-u}h(x,t,\e^us,s\e^u\nabla u +\e^{u}\eta),\\
&\widehat{\bm{m}}=\e^{-u}\bm{m},\widehat{g}(x,t,s)=g(x,t,\e^{-u} s),
\notag\\
&\widehat{\psi}(x,t,s)=\e^{-u}\psi(x,t,\e^us)-Kus,
\notag\\
&\widehat{L}_t[v]= v_t-\divi \ (a\nabla v)+\widehat{\bm{b}}\nabla v +\widehat{c}v,\notag\\
&\widehat{N}[v]=\widehat{h}(x,t,v,\nabla v)+\widehat{\bm{m}}\cdot\nabla (\widehat{g}(x,t,v)).
\end{align*}
Putting $w=\e^uv$ into \eqref{L1}, we obtain \eqref{hat} with $\beta=\e^u$.

Note that 
$-\widehat{h}(x,t,s,\eta)s\leq \mu(x)s^2:=\widehat{\mu}(x)s^2$  holds for all $ (x,t,s,\eta)\in {\Omega}\times \mathbb{R}_{+}\times \mathbb{R}\times \mathbb{R}^n$.
Hence, the structural conditions \textbf{(A2-1)} 
are satisfied.

We verify that \textbf{(A3-2)} is also satisfied. Indeed, it suffices to note that
\begin{align*}
\big(\widehat{\psi}(x,t,s)+s\widehat{\bm{b}}\cdot \bm{\nu}
\big)s
=&\e^{-u}\bigg( \psi(x,t,\e^us)+\bm{n}\cdot \bm{\nu} (\e^us)
\bigg)s-2a(\nabla u\cdot \bm{\nu})s^2\notag\\
\geq &-2Kaus^2\\
\geq &0,\forall (x,t,s)\in {\Omega}\times \mathbb{R}_{+}\times \mathbb{R}.
\end{align*}
Obviously, \textbf{(A1-4)} is satisfied with $\mu_3(x,t)=-2Kau+|\widehat{\bm{b}}|$.

Now note that
\begin{align*}
\widehat{c}- \divi \ \widehat{\textbf{b}}-\widehat{\mu}
=&( -  \divi \ (a\nabla u)-a|\nabla u|^2+\bm{b}\cdot\nabla u+c)-(\divi \ {\textbf{b}} -2\divi \ (a\nabla u))-{\mu}\notag\\
=&-\divi \ {\textbf{b}}+  \divi \ (a\nabla u)-a|\nabla u|^2+\bm{b}\cdot\nabla u+c-{\mu}\notag\\
=& \underline{c}>0,\forall (x,t)\in {\Omega}\times \mathbb{R}_{+}.
\end{align*}
By the result of Theorem \ref{main result} (i), we have
\begin{align*}
\|\e^{-u(\cdot)}w(\cdot,T)\|_{L^1( \Omega)}
\leq  &\e^{-\underline{c}T}\|{\e^{-u}w^0}\|_{L^1( \Omega)}+ \| \e^{-u}{d}\|_{L^1(\partial\Omega\times({0,T}))}+\| \e^{-u}{f}\|_{L^1( Q_T)},
\end{align*}which and $u\leq 0$ result in:
\begin{align*}
\|w(\cdot,T)\|_{L^1( \Omega)}
\leq  &\e^{\|u\|_{L^{\infty}(\Omega)}}\e^{-\underline{c}T}\|{w^0}\|_{L^1( \Omega)}
+ \!\e^{\|u\|_{L^{\infty}(\Omega)}}\| {d}\|_{L^1(\partial\Omega\times({0,T}))}\!+\!\e^{\|u\|_{L^{\infty}(\Omega)}}\| {f}\|_{L^1( Q_T)}.
\end{align*}

\section{Technical lemmas used in Section~\ref{ISS Orlicz}}
\begin{lemma}\cite{Martinez:2008}\label{Lemma 9} The function $\phi^{-1}$ satisfies the
following structural condition:
\begin{align}
\frac{1}{\delta_{1}}\leq \frac{s(\phi^{-1})'(s)}{\phi^{-1}(s)}\leq
\frac{1}{\delta_0},\ \ \forall  s>0.\notag
\end{align}
\end{lemma}

\begin{lemma}\cite[Page 264]{AF2003}\label{Young}
$ab\leq \Phi(a)+\widetilde{\Phi}(b)$ for all $a,b\in \mathbb{R}_{\geq 0}$.
\end{lemma}
\begin{lemma} \label{property of inverse of phi}The following results hold true:
\begin{enumerate}[(i)]
 \item $\frac{\delta_0(1+\delta_1)}{\delta_1(1+\delta_0)} \min\{k^{1+\frac{1}{\delta_0}},k^{1+\frac{1}{\delta_1}}\}\widetilde{\Phi}(s)\leq \widetilde{\Phi}(ks)\leq \frac{\delta_1(1+\delta_0)}{\delta_0(1+\delta_1)} \max\{k^{1+\frac{1}{\delta_0}},k^{1+\frac{1}{\delta_1}}\}\widetilde{\Phi}(s), \ \forall  k,s\geq
 0$;
 \item $ab\leq \varepsilon \Phi(a)+C(\varepsilon)\widetilde{\Phi}(b) $, $\forall a,b\geq 0,\forall \varepsilon\in (0,1)$, where $C(\varepsilon)=\frac{\delta_1(1+\delta_0)}{\delta_0(1+\delta_1)} \bigg(\frac{1+\delta_0}{1+\delta_1} \varepsilon\bigg)^{-\frac{1+\delta_0}{\delta_0(1+\delta_1)}}>0$;
     \item $\widetilde{\Phi}(\phi(s))=s\phi(s)-\Phi(s)\leq \delta_1 \Phi(s)$, $\forall s\geq 0$.
   \end{enumerate}
\end{lemma}
\begin{pf*}{Proof:} {Note first that in }Lemma~\ref{Lemma 9}, $\phi^{-1}$ also satisfies the
Tolksdorf's condition. Applying Lemma \ref{property of phi} (iv), we obtain (i).

Given $\varepsilon\in (0,1)$, let $ \varepsilon_0=\big(\frac{1+\delta_0}{1+\delta_1} \varepsilon\big)^{\frac{1}{1+\delta_1}}\in (0,1)$ and $C(\varepsilon)=\frac{\delta_1(1+\delta_0)}{\delta_0(1+\delta_1)} \big(\frac{1+\delta_0}{1+\delta_1} \varepsilon\big)^{-\frac{1+\delta_0}{\delta_0(1+\delta_1)}}>0$. By Lemma \ref{Young}, for any $a,b\geq 0$, we have
\begin{align*}
 ab=\varepsilon_0a\cdot \frac{b}{\varepsilon_0}
\leq  \frac{1+\delta_1}{1+\delta_0}\varepsilon_0^{1+\delta_1}G(a)+\frac{\delta_1(1+\delta_0)}{\delta_0(1+\delta_1)} \bigg(\frac{1}{\varepsilon_0}\bigg)^{1+\frac{1}{\delta_0}}\widetilde{\Phi}(b)
=\varepsilon G(a)+C(\varepsilon)\widetilde{\Phi}(b),
\end{align*}
which is the result in (ii).

Finally, we refer to \cite{Martinez:2008} for the proof of (iii). 
$\hfill\blacksquare$
\end{pf*}
\begin{lemma} \label{property of inverse of phi'}The following results hold true:
\begin{enumerate}[(i)]
\item for any $u\in L^{\Phi}(\omega)$, it holds:
\begin{align*}
     \min\limits_{i\in\{0,1\}}\bigg\{\bigg(\frac{1+\delta_0}{1+\delta_1}
             \int_{\omega}\Phi(|u|)\diff{y}\bigg)^{\frac{1}{1+\delta_i}}\bigg\}
\leq \| u\|_{L^{\Phi}(\Omega)}
\leq  \max\limits_{i\in\{0,1\}}\bigg\{\bigg(\frac{1+\delta_1}{1+\delta_0}\int_{\omega}\Phi(|u|)\diff{y}\bigg)^{\frac{1}{1+\delta_i}}
\bigg\};
\end{align*}
\item $L^{\widetilde{\Phi}}(\omega)$ endowed with the Luxemburg norm is also a Banach space, which is the dual of $L^{\Phi}(\omega)$;
\item for any $u \in L^{\Phi}(\omega)$ and any $v \in L^{\widetilde{\Phi}}(\omega)$, it holds $\left|\int_{\Omega} uv \diff{x} \right|$ $ \leq 2 \| u\|_{L^{\Phi}(\omega)} \| v\|_{L^{\widetilde{\Phi}}(\omega)}.$
\end{enumerate}
\end{lemma}

\begin{pf*}{Proof:} For $ \int_{\omega}\Phi(|u|)\diff{y}=0$, $u=0$ a.e. in $\omega$. Then the results in Lemma \ref{property of inverse of phi'} (i) follows directly. For $\int_{\omega}\Phi(|u|)\diff{y}\neq 0$, letting
$$ k=\max\limits_{i\in\{0,1\}}\bigg\{\bigg(\frac{1+\delta_1}{1+\delta_0}\int_{\omega}\Phi(|u|)\diff{y}\bigg)^{\frac{1}{1+\delta_i}}
\bigg\},$$
we deduce from Lemma~\ref{property of phi}~(iv) that
\begin{align*}
\int_{\omega}\!\Phi\bigg(\!\frac{|u|}{k}\!\bigg)\diff{y}\leq \! \frac{1+\delta_1}{1+\delta_0}\!\max\limits_{i\in\{0,1\}}\{k^{-(1+\delta_i)}\} \int_{\omega}\!\Phi(|u|)\diff{y}\leq \!1.
\end{align*}
Therefore $\|u\|_{L^{\Phi}(\omega)}\leq k$ and the second inequality in Lemma \ref{property of inverse of phi'} (i) follows.

For any $k>0$, we deduce from $ \int_{\omega}\Phi\big(\frac{|u|}{k}\big)\diff{y}\leq 1$ and Lemma \ref{property of phi} (iv) that
\begin{align*}
1\!\geq\! \int_{\omega}\!\Phi\bigg(\frac{|u|}{k}\bigg)\diff{y}\geq \frac{1+\delta_0}{1+\delta_1}\!\min\limits_{i\in\{0,1\}}\{k^{-(1+\delta_i)}\} \! \int_{\omega}\!\Phi(|u|)\diff{y},
\end{align*}
which implies that
$k\geq \min\limits_{i\in\{0,1\}}\bigg\{\Big(\frac{1+\delta_0}{1+\delta_1}\int_{\omega}\Phi(|u|)\diff{y}\Big)^{\frac{1}{1+\delta_i}}
\bigg\}.$
By the definition of $\|u\|_{L^{\Phi}(\omega)}$, the first inequality in Lemma \ref{property of inverse of phi'} (i) follows.

Finally, we refer to \cite[Page 269, 373]{AF2003} for the proof of Lemma~\ref{property of inverse of phi'}~(ii) and (iii).
$\hfill\blacksquare$
\end{pf*}

\end{document}